\documentclass[10pt]{article}

\usepackage{a4,amsmath,amssymb,amscd,graphicx,epsfig,xypic}
\usepackage{vmargin}

\everymath{\displaystyle}

\newtheorem{teo}{Theorem}[section]
\newtheorem{lem}[teo]{Lemma}
\newtheorem{cor}[teo]{Corollary}

\newtheorem{prop}[teo]{Proposition} 
\newtheorem{lem-defi}[teo]{Lemma-Definition}
 \newtheorem{defi}[teo]{Definition} 

\newtheorem{conge}[teo]{Conjecture}
\newtheorem{remark}[teo]{Remark}
\newtheorem{proble}[teo]{Problem}

\newcommand{\mr}{\mathbb{R}}
\newcommand{\mc}{\mathbb{C}}
\newcommand{\mz}{\mathbb{Z}}

\newcommand{\mn}{\mathbb{N}}

\title{Quantum Hyperbolic State Sum Invariants of 3-Manifolds}

\author{St\'ephane Baseilhac and Riccardo Benedetti}

\date{}

\begin{document}

\maketitle

\vspace{1cm}

\noindent Laboratoire E. Picard, CNRS UMR 5580, Universit\'e Toulouse III, F-31062 TOULOUSE

\noindent Email: baseilha@picard.ups-tlse.fr

\vspace{0.3cm}

\noindent Dipartimento di Matematica, Universit\`a di Pisa, Via F. Buonarroti, 2, I-56127 PISA

\noindent Email: benedett@dm.unipi.it

\vspace{1cm}

\begin{abstract} 
\noindent Any triple $(W,L,\rho)$, where $W$ is a compact closed
oriented 3-manifold, $L$ is a link in $W$ and $\rho$ is a flat principal
$B$-bundle over $W$ ($B$ is the Borel subgroup of upper triangular
matrices of $SL(2,\mc)$), can be encoded by suitable 
{\it distinguished} and {\it decorated}
triangulations ${\cal T}=(T,H,{\cal D})$. For each $\cal T$, for each
odd integer $N\geq 3$, one defines a state sum $K_N({\cal T})$, based on 
the
Faddeev-Kashaev quantum dilogarithm at $\omega =\exp( 2\pi i/N)$, such that
$K_N(W,L,\rho)=K_N({\cal T})$ is a well-defined complex valued invariant.
The purely topological, conjectural invariants $K_N(W,L)$ proposed earlier
by Kashaev
correspond to the special case of the {\it trivial} flat bundle.
Moreover, we extend the definition of these invariants to the case of
flat bundles on $W\setminus L$ with not necessarily trivial holonomy
along the meridians of the link's components, and also to $3$-manifolds 
endowed with a $B$-flat bundle and with \emph{arbitrary} non-spherical 
parametrized boundary components. As a matter of fact the distinguished 
and decorated triangulations are strongly reminiscent of 
the way one represents the classical {\it refined scissors congruence} 
class $\widehat{\beta} ({\cal F})$, belonging to the extended Bloch group, 
of any given finite volume hyperbolic 3-manifold $\cal F$ by 
using any hyperbolic ideal triangulation of $\cal F$. We point out
some remarkable specializations of the invariants; among these,
the so called  {\it Seifert-type} invariants, when $W=S^3$: these seem to be 
good candidates in order to fully reconstruct the Jones polynomials in the 
main stream of quantum hyperbolic invariants. 
Finally, we try to set our results against the heuristic backgroud of the 
Euclidean analytic continuation 
of (2+1) quantum gravity with negative cosmological constant, 
regarded as a gauge theory with the {\it non-compact} group $SO(3,1)$ 
as gauge group. 
\end{abstract}

\bigskip

\noindent \emph{Keywords: quantum dilogarithm, hyperbolic 3-manifolds, state sum invariants.}

\section{Introduction}

In a series of papers \cite{1,2,3}, Kashaev
proposed a conjectural infinite family $\{ K_N \}$, $N$ being an odd
positive integer, of complex valued topological invariants of
links $L$ in any oriented compact closed 3-manifold $W$, based on the
theory of quantum dilogarithms at an $N$'th-root of unity 
$\omega = \exp(2\pi i/N)$ \cite{4,5}. These invariants should be 
computed, in a purely
3-dimensional context, as a {\it state sum} $K_N({\cal T})$, supported
by an arbitrary distinguished and decorated triangulation of $(W,L)$,
${\cal T}=(T,H,{\cal D})$ say, if any. ``Distinguished'' means that
$L$ is triangulated by a Hamiltonian subcomplex $H$ of $T$; the rather
complicated decoration $\cal D$ shall be specified later. The main
ingredients of the state sum are the quantum-dilogarithm $6j$-symbols
suitably associated to the tetrahedra of $T$. On the ``quantum'' side,
the state sum is similar to the Turaev-Viro one \cite{6,7}.
On the ``classical'' side, it relates to the computation of the
volume of a hyperbolic 3-manifold by the {\it sum} of the volumes of the
ideal tetrahedra of any of its ideal triangulations.

\noindent Beside a somewhat neglected {\it existence} problem of
such distinguished and decorated triangulations of $(W,L)$, a 
main question left unsettled is just the {\it invariance} of
$K_N({\cal T})$ when $\cal T$ varies. On the other hand, Kashaev
proved the invariance of $K_N({\cal T})$ under certain ``moves'' on
distinguished and decorated triangulations, which gives some evidence
for the conjectured topological invariance.
\medskip

\noindent In \cite{2} he also defined a family $Q_N$ of 
topological invariants for links $L$ in $S^3$ using the solutions of
the Yang-Baxter equation (i.e. the $R-matrices$) derived from the 
pentagon identity satisfied by the quantum-dilogarithm 6j-symbols
\cite{5}. Then he argued that these invariants coincide with the
previous ones in the special case when $W=S^3$. More recently,
Murakami-Murakami \cite{8} have shown that $Q_N$ actually
equals a specific  coloured Jones invariant, getting, by
the way, another proof that it is  a well-defined
invariant for links in $S^3$.

\medskip

\noindent Having reformulated Kashaev's invariants (of links in $S^3$)
within the main stream of Jones polynomials has been an important
achievement, but it also has the negative consequence of
putting aside the original purely 3-dimensional set-up (for links in
an arbitrary $W$), willingly forgetting the complicated and somewhat
mysterious decorations.

\noindent On the contrary, in our opinion, Kashaev's 3-dimensional
set-up deserved to be deeper understood and developed.  As a
by-product, we shall see in \S \ref{GR} that there are good reasons to
believe that it could contain a consistent part of an ``exact
solution'' of the Euclidean analytic continuation of (2+1) quantum
gravity with negative cosmological constant, regarded as a gauge
theory with the {\it non-compact} gauge group $SO(3,1)$ and an action
of Chern-Simons type.  In fact, the so called {\it Volume Conjecture}
\cite{9} about the asymptotic behaviour of $Q_N$ when $N\to \infty$
perfectly agrees, for what concerns the ``real part'', with the
expected ``classical limit'' of this theory \cite[p.77]{10}. In
particular, recall that hyperbolic 3-manifolds are the pure-gravity
(i.e. empty) Euclidean ``classical'' solutions.
\medskip

\noindent The initial aim of the present paper was to work out a
proof that the original 3-dimensional Kashaev set-up actually
produces topological invariants. In fact, after having made it more
``flexible'' (as is demanded by technical and also conceptual
reasons), and having better understood 
the decoration, we have finally obtained the following more
general results that we state here in a somewhat informal way.

\smallskip

\noindent {\bf Notations.} We denote by $W$ a compact closed
oriented 3-manifold, $L$ is a link in $W$, $U(L)$ is a tubular
neighbourhood of $L$ in $W$, and $M= W\setminus \ {\rm Int}\ U(L)$. 
So $M$
has $n$ toral boundary components, where $n$ is the number of components
of $L$; we denote them by $L_i$, $i=1,\dots,n$.
Let $Z$ be any compact oriented 3-manifold with non-empty boundary
$\partial Z = \coprod_{i=1}^n T_i^2$ consisting of $n$ toral components 
$T_i^2$, and  $Y$ any compact oriented 3-manifold with non-empty 
boundary $\Sigma$, so that each of its connected components $\Sigma_i$ 
has genus $g(\Sigma_i)\geq 1$. 

\noindent Given $Z$, $S=\{s_i\}$ denotes a system of 
{\it essential}, that is non contractible, simple closed curves, one on 
each $T_i^2$. \hfill\break
\noindent $B=B(2,\mc)$ is the Borel subgroup of $SL(2,\mc)$ of upper 
triangular matrices, and $\rho$ denotes either an equivalence class of 
{\it flat} principal $B$-bundles on $Z$ (resp. $Y$), i.e. a bundle 
endowed with a flat connection, or, equivalently, an element of 
$\chi_B(Z) = Hom(\pi_1(Z),B)/B$ (resp. $\chi_B(Y) = Hom(\pi_1(Y),B)/B$), 
$B$ acting by inner automorphisms. One usually calls $\chi_B(Y)$
the \emph{character variety of $Y$ for $B$}.

\noindent Finally, we denote by $\phi$ a \emph{symplectic} 
parametrization 
of $\Sigma$ (see Section \ref{ZS} for the precise definition).   
The triples $(Z,S,\rho)$ and $(Y,\phi,\rho)$ are regarded up to 
the natural equivalence relation induced by orientation preserving 
homeomorphisms.
\medskip

\noindent  We first consider triples $(Z,S,\rho)$.

\begin{teo}\label{t1} A family of complex valued 
invariants $\{K_N(Z,S,\rho)\}$, $N$ being any odd positive integer, is defined. The triple $(Z,S,\rho)$ can be encoded
by suitably decorated ``distinguished'' triangulations 
$(T,H,{\cal D})$, and $K_N(Z,S,\rho)$ can be computed as a
state sum $K_N(T,H,{\cal D})$ whose value is invariant when 
$(T,H,{\cal D})$ varies. The original Kashaev's purely topological 
invariants $K_N(W,L)$ correspond to the case when $Z=M$, each $s_i$ is a 
meridian of $L_i$
and $\rho$ equals the class of the trivial flat bundle, or,
equivalently, equals the constant representation $\rho=1$.
\end{teo}

\begin{remark} \label{r1} {\rm There are two ways to look at $(Z,S,\rho)$, 
whence at these invariants.

(1) Let $W$ be the closed manifold obtained by {\it Dehn filling}
of each boundary component of $Z$, along the curves of $S$. If $L$
denotes the link of surgery ``cores'', then $S$ becomes a family of
meridians of the $L_i$'s, and $\rho$ is a flat bundle on $W\setminus L$
with, in general, a {\it non-trivial holonomy} along these meridians.
In classical $(2+1)$ gravity an analogous situation arises in presence of 
closed 
lines of massive particles of a universe (see for instance \cite{11}).

(2) We are actually considering $Z$ with {\it parametrized} boundary.
More precisely, fix a {\it base} solid torus ${\cal H}_1$ with meridian
$m$.  Then each boundary component $T_i$ of $Z$ is parametrized by an
oriented diffeomorphism $g_i:\partial {\cal H}_1 \to T_i$ such that
$g_i(m)=s_i$;  $g_i$ is well-defined up to automorphisms of $\partial
{\cal H}_1$ which extend to the whole ${\cal H}_1$.}
\end{remark}

\noindent The symplectic parametrizations of surfaces cited above are 
defined in the spirit of Remark \ref{r1} (2), starting with {\it base}
handlebodies of genus $g(\Sigma_i)$ with fixed systems of meridians
and so on. Then we also have:

\begin{teo}\label{t2} One can define a family of complex valued 
invariants $\{K_N(Y,\phi,\rho)\}$, $N$ being any odd positive integer, 
which 
specializes to the previous one when $\partial Y$ is a union of tori.
\end{teo}
\noindent Our final set-up shall be slightly different from
the original Kashaev one. So we prefer to reach it step by step.
\smallskip

\noindent 
In \S \ref{state-sum1}-\ref{state-sum2}, we 
shall prove Theorem \ref{t1} in a particular case: $Z=M$, each $s_i$ is a 
meridian of $L_i$ and $\rho$ is defined on the whole
of $W$, i.e. $\rho$ is trivial along the $s_i$'s; that is we actually 
construct invariants $K_N(W,L,\rho)$.  Note that
the original Kashaev invariants already are specializations of them 
(with $\rho$ trivial). To treat this case we shall need a milder 
modification of the original Kashaev set-up. The related {\it existence 
problem} of distinguished and decorated 
triangulations is solved in  \S \ref{dist}-\ref{deco}.  
Some basic properties of these invariants are settled in 
\S \ref{properties}: ``projective invariance'', and ``duality'' (which in 
particular describes the behaviour of the invariants when the orientation 
of $W$ changes).

\smallskip

\noindent  Next, we treat the case when $\rho$ does not necessarily 
extend to the whole of $W$.
This is done in \S \ref{ZS} through a trick, which consists roughly 
in cutting up pieces of a distinguished triangulation of $(W,L)$, turning 
it into a decorated cell decomposition of $Z$ with tetrahedron-like
building blocks, which may be used to again define our previous state sums. 
Thanks to our definition of parametrized surface $(\Sigma,\phi)$, there is 
a natural way to 
recover the set-up used in the proof of Theorem \ref{t1} from triples 
$(Y,\phi,\rho)$, thus obtaining Theorem \ref{t2}.

\smallskip

\noindent In \S \ref{specialize}, we consider some remarkable
specializations of the invariants $K_N(Z,S,\rho)$; among these 
there are the so called \emph{Seifert-type} invariants of links in 
$\mathbb{Z}$-homology spheres. For $S^3$, they seem to be good candidates
in order to fully recover the coloured Jones polynomials in terms of the 
quantum hyperbolic state sum invariants (generalizing \cite{8}), and to 
give them a new interpretation.
One finds a preliminary discussion at the end of \S \ref{specialize},
our work on this matter being in progress (cf. \cite{12}).   

\smallskip

\noindent We stress a qualitative contribution of the present paper 
in understanding the nature of the distinguished and decorated 
triangulations: it turns out that they are strongly 
reminiscent (up to a sort of reduction mod $N$) of the way one represents 
the classical {\it refined scissors congruence}
class $\widehat{\beta} (\mathcal {F})$ (belonging to the extended Bloch group
$\widehat{B} (\mc)$) of any given finite volume hyperbolic 3-manifold 
$\mathcal{F}$, by using any of its hyperbolic ideal triagulations 
\cite{13,14}.
It is known that the classical Rogers dilogarithm lifts, with its 
fundamental {\it five-terms
identity}, to an analytic function $R$ defined on $\widehat{B} (\mc)$, and
that
$$R(\widehat{\beta} ({\cal F}))= i 
\left( Vol({\cal F})+iCS({\cal F}) \right),$$
\noindent  where $CS$ denotes the Chern-Simons invariant \cite{14,15}.
Moreover one knows that the $N$-dimensional quantum dilogarithm 6j-symbols 
and their pentagonal 
relations recover respectively, in the limit $N\to \infty$, the classical 
Euler dilogarithm and the Rogers five-terms identity (see \cite{4,5} and, for a detailed account, \cite{16}). All this should indicate
the right conceptual framework supporting the 
so called {\it volume conjecture} and natural generalizations, also 
involving the Chern-Simons invariant; roughly speaking one expects that: 
\medskip

\noindent {\bf Qualitative Hyperbolic  $VCS$-Conjecture.} 
{\it When $(W,L)$ admits some ``hyperbolization'' $\mathcal{F}$ 
{\rm (for instance: 
$\mathcal{F}= W\setminus L$ is a
complete hyperbolic manifold or $\mathcal {F}=W$ is hyperbolic and $L$ is
a geodesic link, as it arises from hyperbolic Dehn surgery)}, then one 
can recover 
$VCS(\mathcal {F})=Vol(\mathcal {F})+iCS(\mathcal {F})$ by means of the
asymptotic behaviour for $N\to \infty$ of the invariants $K_N(W,L,\rho)$.}

\medskip 

\noindent Nevertheless, some substantial facts are not yet well understood, 
and this matter needs further investigations (see also \S \ref{charges}, 
\S \ref{specialize} and \S \ref{GR}). 
\smallskip

\noindent Finally in \S \ref{GR} we try to set our results against the
heuristic backgroud of the Euclidean analytic continuation of
(2+1) quantum gravity with negative cosmological constant, regarded as
a gauge theory with the {\it non-compact} group $SO(3,1)$ as gauge
group, as defined by Witten \cite{10}. We shall develop a few
speculations about a path integral interpretation of our state sums and
their relation with the absolute torsions of Farber-Turaev \cite{17}.
\smallskip

\noindent The aim of \S \ref{App} is to provide the reader a
precise statement of all the facts we use concerning our ``quantum
data'', allowing a substantially self-contained reading.  We shall
mostly insist on the geometric interpretation of the basic properties,
showing how the decorations we use in this paper do arise. All the
statements given without proof in the Appendix are carefully treated
in \cite{16}.

\section{Distinguished singular triangulations of $(W,L)$} \label{dist}

Let us recall some facts about standard spines of 3-manifolds
and their dual {\it ideal} triangulations \cite{18,19}. A reference is \cite{20}, but
one finds also a clear discussion about this material in \cite{6} 
(note that sometimes the terminologies do not agree). Using the 
Hauptvermutung, we
will freely intermingle the differentiable, piecewise linear and 
topological
viewpoints for $3$-dimensional manifolds. 
\medskip

\noindent Consider a tetrahedron $\Delta$ and let $C$ be the interior
of the 2-skeleton of the dual cell-decomposition. 
A {\it simple} polyhedron $P$ is a finite 2-dimensional 
polyhedron such that each point of $P$ has a neighbourhood which can be
embedded into $C$. A simple polyhedron is {\it standard} (in \cite{6}
one uses the term {\it cellular}) if all the
components of the natural stratification of $P$ given by singularity are 
open cells. Depending on the dimension, we call these components 
{\it vertices, edges} and {\it regions} of $P$. 

\noindent Every compact 3-manifold
$Y$ (which for simplicity we assume connected) with non-empty boundary
$\partial Y$ has {\it standard spines}, that is standard polyhedra $P$
embedded in Int $Y$ such that $Y$ collapses onto $P$ (i.e. $Y$ is a
regular neighbourhood of $P$). Standard spines of oriented 3-manifolds
are characterized among standard polyhedra by the property of carrying
an {\it orientation}, that is a suitable ``screw-orientation'' along
the edges \cite[Def. 2.1.1]{20}.  Such an oriented 3-manifold $Y$ can
be reconstructed (up to orientation preserving homeomorphisms) from
any of its oriented standard spines. From now on we assume that $Y$ is
oriented.

\medskip

\noindent A {\it singular} triangulation of a polyhedron $Q$ is a
triangulation in a weak sense, namely self-adjacencies and multiple
adjacencies are allowed. For any $Y$ as above, let us denote by $Q(Y)$
the polyhedron obtained by collapsing each component of $\partial Y$
to a point. An {\it ideal triangulation} of $Y$ is a singular 
triangulation
$T$ of $Q(Y)$ such that the vertices of $T$ are precisely the points of 
$Q(Y)$ which correspond to the components of $\partial Y$.

\medskip

\noindent For any ideal triangulation $T$ of $Y$, the 2-skeleton
of the  \emph{dual} cell-decomposition of $Q(Y)$ is a standard spine
$P(T)$ of $Y$. This procedure can be reversed, so that we can associate
to each standard spine $P$ of $Y$ an ideal triangulation $T(P)$ of $Y$
such that $P(T(P))=P$. Thus standard spines and ideal triangulations
are dual equivalent viewpoints which we will freely intermingle.
Note that, by removing small neigbourhoods of the vertices of $Q(Y)$,
any ideal triangulation leads to a cell-decomposition of $Y$ by 
{\it truncated tetrahedra} which induces a singular triangulation of the
boundary of $Y$.

\begin {remark} 
\label{rem1} {\rm For the following facts, a reference is 
\cite[Ch. E]{21}. The name ``ideal triangulation'' is 
inspired by the geometric triangulations of a non compact hyperbolic
3-manifold with finite volume (for example, the complement of a
hyperbolic link in $S^3$) by ideal hyperbolic tetrahedra (possibly
partially flat). It is a standard result of Epstein-Penner
\cite{22} that such triangulations exist. A geometric ideal 
triangulation can be regarded as a special topological ideal
triangulation admitting suitable decorations by the {\it moduli} of
the corresponding hyperbolic ideal tetrahedra.  Then, the volume of
the 3-manifold can be expressed as the sum of the volumes of the ideal
tetrahedra of any of its geometric triangulations, and it can be
computed in terms of the moduli via the Bloch-Wigner dilogarithm
function \cite{23}. Of course the volume, which by Mostow's rigidity
theorem is a topological invariant, does not depend on the specific
geometric ideal triangulation used to compute it. In a sense, the
quantum hyperbolic state sums we are concerned with can be considered
as ``quantum'' deformations of this ``classical'' hyperbolic
situation, which make sense for arbitrary 3-manifolds.}
\end{remark}

\noindent Consider now our closed 3-manifold $W$. For any $r\geq 1$
let $W'_r = W\setminus rD^3$, that is the manifold with $r$ spherical
boundary components obtained by removing $r$ disjoint open balls from $W$.
Clearly $Q(W'_r)=W$ and any ideal triangulation of $W'_r$ is a singular
triangulation of $W$; moreover all singular triangulations of $W$ 
are obtained in this way. We shall adopt the following terminology.

\begin{defi}{\rm A singular triangulation of $W$ is simply called
a {\it triangulation}. Ordinary triangulations (where neither 
self-adjacencies
nor multi-adjacencies are allowed) are said to be {\it regular}.
An {\it almost-regular} triangulation of $W$ is a triangulation 
having the same vertices as a regular one. This is equivalent to saying
that there exists a regular triangulation of $W$ with the same number
of vertices. Given a triangulation $T$, $r_i=r_i(T)$, $i=0,1,2,3$,
shall denote the number  of vertices, edges, faces, tetrahedra of $T$. }
\end{defi} 

\noindent The main advantage in using singular
triangulations (standard spines) instead of only ordinary triangulations 
consists of the
fact that there exists a {\it finite} set of moves which are sufficient
in order to connect (by means of finite sequences of these moves) 
singular triangulations (standard spines) of the same manifold.
On the contrary, if we pretend to keep ordinary 
triangulations at each steps, we are forced to consider 
an {\it infinite} set of moves (see \cite{6}). 

\medskip

\begin{figure}[ht]

\begin{center}

\scalebox{1}{\input{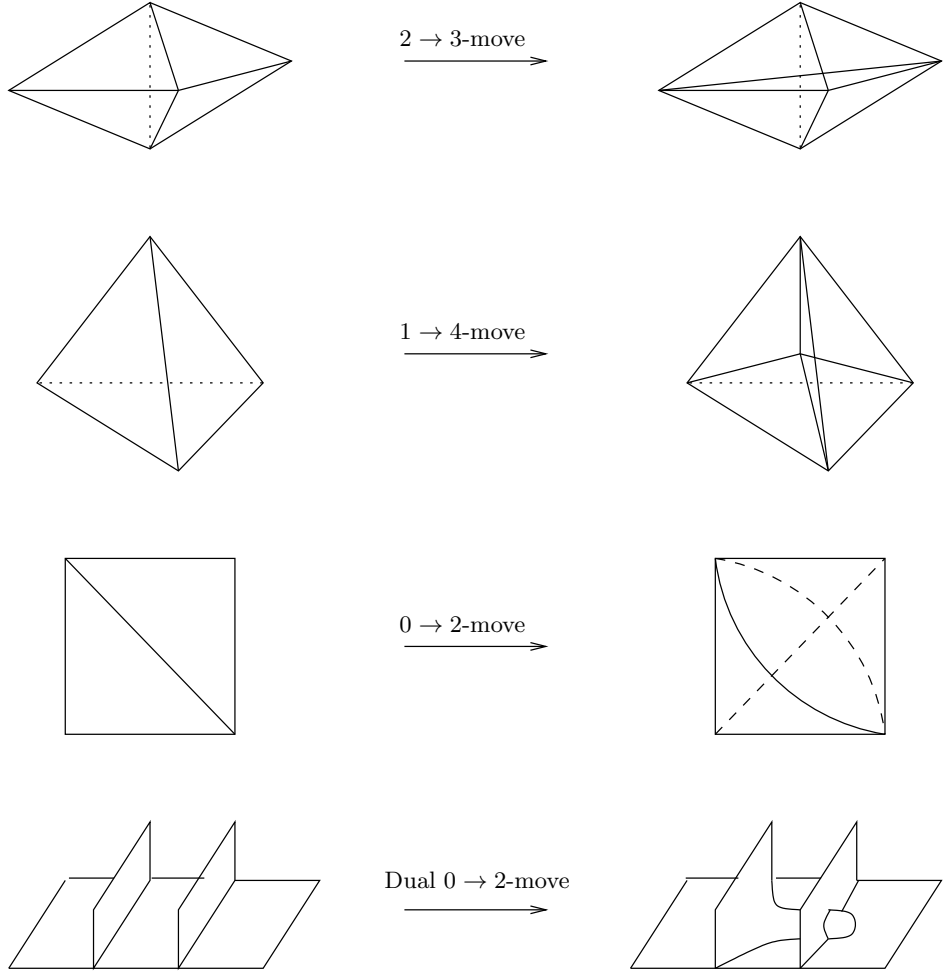}}

\end{center}

\caption{\label{fig1} the moves on triangulations, 
and the dual $0 \to 2$ move on simple spines}

\end{figure}

\noindent Let us recall two elementary moves on triangulations (spines)
that we shall use throughout the paper; see Fig. \ref{fig1}.

\smallskip

\noindent {\bf The $2\to 3$ move.} Replace the 
triangulation $T$ of a portion of $Q(Y)$ made by the union of $2$ 
tetrahedra
with a common 2-face $f$ by the triangulation made by $3$ tetrahedra
with a new common edge which connect the two vertices opposite to $f$.

\smallskip

\noindent {\bf The $1\to 4$ move.} Add  a new vertex
in the interior of a tetrahedron $\Delta$ of $T$ and make from it the cone
over the triangulated boundary of $\Delta$. The dual spine $P'$ of the
triangulation $T'$ thus obtained is a spine of $Y\setminus D^3$, 
where $D^3$ is an open ball in the interior of $Y$.

\medskip

\noindent The $2\to 3$ and $3\to 2$ moves 
can be easily reformulated in dual terms 
(see for instance \cite[p.15]{20}, where they are called ``MP-moves''; 
in Fig. \ref{fig:slide}, we show their branched versions). 
Standard spines of the same $Y$ with at least two
vertices (which, of course, is a painless requirement) may always be
connected by 
the (dual) move $2\to 3$ and its inverse. In particular, any almost-regular
triangulation of $W$ can be obtained from a regular
one via a finite sequence of $2\to 3$ or $3\to 2$ moves.   

\noindent In order to handle 
singular triangulations of a closed manifold $W$, we also need a move which
allows us to vary the number of vertices. Although this is not the
shortest way (the so-called {\it bubble} move makes a hole in any $Y$
by introducing only two more vertices, see Proposition \ref{vertex}), 
the $1\to 4$ move (and its inverse)
shall be convenient for our purposes. Let us describe the dual 
spine $P'$ of  $Y\setminus D^3$. Consider the vertex $v_0$ of $P$ 
(which is dual
to $T$) contained in $\Delta$. We are removing a ball around of $v_0$.
The boundary $S^2(v_0)$ of this ball is a portion of $P'$ and the natural 
cell 
decomposition of $P'$ induces on it a cell decomposition which is 
isomorphic
to the ``bi-dual'' tetrahedron of $v_0$; that is it is isomorphic to the
cell decomposition of $\partial \Delta$ dual to the standard triangulation. 

\noindent For technical reasons we shall need a further move which we
recall in Fig. \ref{fig1} in terms of standard spines and that we
denote by $0\to 2$ move. It is also known as {\it lune} move and is
somewhat similar to the second Reidemeister move on link diagrams. Note
that the inverse of the lune move is not always admissible because one
could lose the standardness property when using it.  However we shall
need only the ``positive'' lune move.

\noindent The following technical result due to Makovetskii \cite{24} 
shall be necessary.

\begin{prop} \label{chemin} Let $P$ and $P'$ be standard spines of $Y$. Then
there exists a spine $P''$ of $Y$ such that $P''$
can be obtained from both $P$ and $P'$ via a finite sequence of $0 \to 2$ 
and $2 \to 3$ moves (we stress that they are all positive moves).
\end{prop}

\begin{defi} {\rm A {\it distinguished} triangulation $(T,H)$ of $(W,L)$ has 
the property that $H$ is 
Hamiltonian, that is at each vertex of $T$ there are exactly two ``germs''
of edges of $H$. Of course, the two germs could belong to the same edge
of $H$.}
\end{defi}

\noindent It is convenient to give a slightly different
description of the distinguished singular triangulations in terms of
spines. 

\begin{defi} \label{adapte} {\rm Let $Y$ be as before. 
Let $S$ be any finite family of $r$ 
disjoint simple closed curves on $\partial Y$. We say that $Q$ is a
{\it quasi-standard} spine of $Y$ {\it relative} to $S$ if:

(i) $Q$ is a simple polyhedron with boundary $\partial Q$ consisting
of $r$ circles. These circles bound (unilaterally) $r$ annular regions
of $Q$. The other regions are cells.

(ii) $(Q,\partial Q)$ is properly embedded in $(Y,\partial Y)$ and 
transversely intersects $\partial Y$ at $S$.

(iii) $Q$ is is a spine of $Y$.

\noindent Note that $\widetilde{P} = Q\setminus \{$annular regions$\}$ 
is a {\it simple} spine of $Y$.}  
\end{defi}

\begin{lem} \label{existadapte} Let $Y$ and $S$ be as before. 
Quasi-standard spines of $Y$
relative to $S$ do exist.
\end{lem}

\noindent {\it Proof.} Let $\widetilde{P}$ be any standard spine of $Y$. 
Consider 
a {\it normal} retraction $r:Y\to \widetilde{P}$. Recall that $Y$ is
the mapping cilynder of $r$; for each region $R$ of $\widetilde{P}$,
$r^{-1}(R) = R\times I$; for each edge $e$, $r^{-1}(e) = e\times \{$a
``tripode''$\}$; for each vertev $v$, $r^{-1}(v) = \{$a
``quadripode''$\}$. We can assume that $S$ is in ``general position''
with respect to $r$, so that the mapping cylinder of the restriction
of $r$ to $S$ is a simple spine of $Y$ relative to $S$ (with the
obvious meaning of the words); possibly after doing some $0\to 2$
moves we then obtain a quasi-standard Q.

\begin{defi} \label{type} {\rm Consider $Y=M$ and $S$ formed by the union 
of $r_i\geq 1$
parallel copies of the meridian $m_i$ of the component $L_i$ of $L$,
$i=1,\dots,n$. A spine of $M$ {\it adapted} to $L$ of {\it type}
$(r_1,\dots,r_n)$ is a quasi-standard spine of $M$ relative to such an $S$.}
\end{defi}

\noindent It is clear that Proposition \ref{chemin} extends to the 
``adapted'' setting. In particular, we have:

\begin{prop} \label{chemintype} Let $P$ and $P'$ be quasi-standard spines 
of $Y$ relative to $S$ and of type $r=(r_1,\ldots,r_n)$. Then there exists 
a spine $P''$ of $Y$ relative to $S$ and of type $r$ such that $P''$
can be obtained from both $P$ and $P'$ via a finite sequence of 
$0 \to 2$ and $2 \to 3$ moves, and at
each step we still have spines of $Y$ adapted to $S$ and of type $r$.
\end{prop}

\begin{remark} \label{trig} {\rm If $Q$ is a spine of $M$ adapted to $L$ as before,
then by filling each boundary component of $Q$ by a 2-disk, we get a
standard spine $P=P(Q)$ of $W'_r$, $r=\sum _i r_i$, and the dual
triangulation $T(P)$ of $W$ naturally contains $L$ as a
Hamiltonian subcomplex; that is we have obtained a distinguished
triangulation of $(W,L)$.  Vice-versa, starting from any
$(T,H)$, by removing an open disk in the dual region
to each edge of $H$, we pass from $P=P(T)$ to a spine $Q=Q(P)$ of $M$
adapted to $L$, of some type. So they are equivalent viewpoints.}
\end{remark}

\noindent As an immediate corollary we have:

\begin{cor} \label{exist1} For any $r\geq n$ there exist 
distinguished 
triangulations of $(W,L)$ with $r$ vertices.
\end{cor}

\noindent Let us analyze the moves on distinguished triangulations.
The $2\to 3$ and $1\to 4$ moves specialize
to moves $(T,H)\to (T',H')$ between distinguished triangulations as
follows.

\noindent Let us denote  by $T\to T'$ any such a move. If $(T,H)$ is
a distinguished triangulation of $(W,L)$, we want to get moves
$(T,H)\to (T',H')$. In the $2\to 3$ case there is nothing to do because
$H'=H$ is still Hamiltonian. In the $1\to 4$ case, we assume that an edge
$e$ of $H$ lies in the boundary of the involved tetrahedron $\Delta$; $e$ 
lies in the
boundary of a unique 2-face $f'$ of $T'$ containing the new
vertex. Then we get the Hamiltonian $H'$ just by replacing $e$ by the
other two edges of $f'$. We have also: 

\begin{lem} Let $(T,H)$ be a distinguished triangulation of
$(W,L)$ and $T\to T'$ a $0\to 2$ move. Then it can be completed to
a move $(T,H)\to (T',H')$.
\end{lem}

\noindent {\it Proof.} This is clear if we think in dual terms. 
If $Q$ is a spine
adapted to $L$ and $Q\to Q'$ is the $0\to 2$ move, the dual regions in
$P(Q)$ to the edges in $H$ ``persist'' in $P(Q')$ so we find $H'$.

\medskip

\noindent Finally we can solve {\it half} of the {\it existence} 
problem mentioned in the introduction.

\begin{prop} There exist almost-regular distinguished 
triangulations of $(W,L)$.
\end{prop}

\noindent {\it Proof.} 
It is enough to remark that any distinguished triangulation,
which exists by Corollary \ref{exist1}, can be made almost-regular
after a finite number of $1\to 4$ moves.

\section{Decorations}\label{deco}

In this section we have to properly define, and possibly better
understand, the {\it decorations} $\mathcal{D}=(b,z,c)$ of 
distinguished triangulations of $(W,L)$.

\subsection{Branchings} \label{branchings}

\noindent Let $P$ be a standard spine of $Y$ and consider as usual
the dual ideal triangulation $T=T(P)$. A {\it branching} $b$ of $T$ is a
system of orientations on the edges of $T$ such that each ``abstract'' 
tetrahedron
of $T$ has one source and one sink on its 1-skeleton. This is
equivalent to saying that, for any 2-face $f$ of $T$, the
edge-orientations do not induce an orientation of the boundary of
$f$.

\medskip 

\noindent In dual terms, a branching is a system of orientations on the
regions of $P$ such that for each edge of $P$ we have the same induced orientation only twice. In particular, note that each edge of
$P$ has an induced orientation. 
\smallskip

\noindent In the original set-up of \cite{1} one used {\it regular}
triangulations $T$ of $W$, with a given total ordering on the set of
vertices; in fact the role of the ordering is just to define a branching 
via the natural
rule: ``on each edge, go from the smaller vertex towards the bigger
one''. Note that one could not exclude, a priori, that after some 
negative moves, one eventually reaches veritable singular triangulations,
for which such a total ordering no longer induces a branching.
\medskip

\noindent Branchings, mostly in terms of spines,
have been widely studied in \cite{20} (see also \cite{25}). 
They are rich structures: 
a branching of
$P$ allows us to give the spine the extra structure of an embedded and
oriented (hence normally oriented) {\it branched surface} in Int($Y$); 
by the way, this also justifies the name. Moreover a branched $P$ 
carries a suitable positively transverse {\it combing} of $Y$. 

\noindent  We recall here part of their combinatorial content. 
A branching $b$ allows to define an orientation on \emph{any} cell of $T$, 
not only on the edges. Indeed, consider any ``abstract'' tetrahedron 
$\Delta$. For each vertex of $\Delta$ consider the number of incoming
$b$-oriented edges in the 1-skeleton. This gives us an 
ordering $b_{\Delta}:\{0,1,2,3\}\to V(\Delta)$ of the vertices which 
reproduces the
branching on $\Delta$, according to the former rule.
This gives us a {\it base} vertex $v_0=b_{\Delta}(0)$
and an ordered triple of edges emanating from $v_0$, whence an orientation 
of $\Delta$.
Note that this orientation may or may not agree with the orientation of $Y$; 
in the first case
we say that $\Delta$ is of index $-1$, and it is of index $1$ otherwise.    

\noindent To orient 2-faces we work in a similar way on the boundary of 
each ``abstract'' 2-face $f$.
We get an ordering $b_f:\{0,1,2\}\to V(f)$, a {\it  base} vertex
$v_0=b_f(0)$, and finally an orientation of $f$. This 2-face orientation
can be described in another equivalent way. Let us consider the 1-cochain 
$s_b$ such that $s_b(e)=1$ for each $b$-oriented edge. Then there is a 
unique way to orient any 2-face $f$ such that the coboundary 
$\delta s_b(f)=1$. 

\noindent The corresponding dual orientation on the edges of $P$ is
just the induced orientation mentioned in \S \ref{dist}. 

\medskip

\noindent {\bf Branching's existence and transit.} This matter is 
carefully 
analyzed in \cite{20}. In general, 
a given ideal triangulation of $Y$ could admit no branching, 
but there exist branched ideal triangulations of any $Y$. More precisely, 
given {\it any} 
system
of edge-orientations $g$ on $T$ and any move $T \to T'$, a \emph{transit}
$(T,g)\to (T',g')$ is given by any system $g'$ of edge-orientations on
$T'$ \emph{which agrees with $g$ on the ``common'' edges}. 
In \cite[Th. 3.4.9]{20} one proves:

\begin{prop} \label{brancheable} For any $(T,g)$ there exists a 
finite sequence of  $2\to 3$
transits such that the final $(T',g')$ is actually branched.
\end {prop}   

\begin{figure}
  \centerline{\psfig{figure=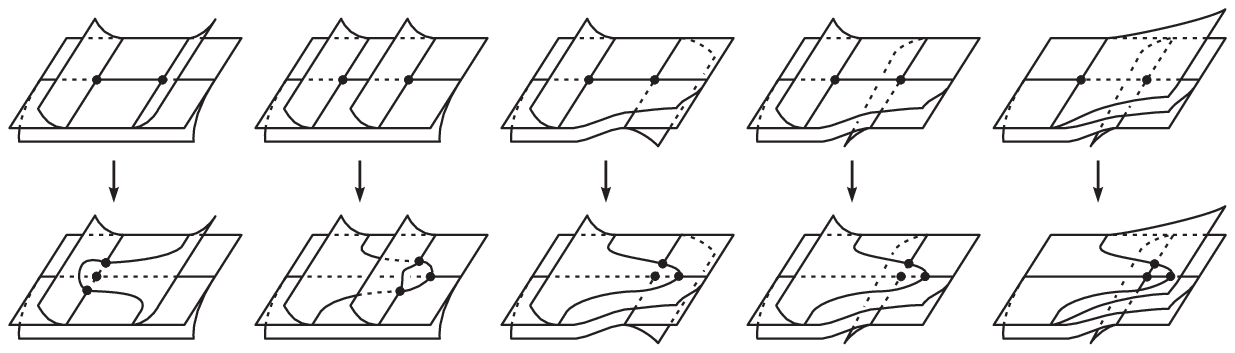}}
  \caption{$2\to 3$ sliding moves.}\label{fig:slide}
\end{figure}

\begin{figure}
  \centerline{\psfig{figure=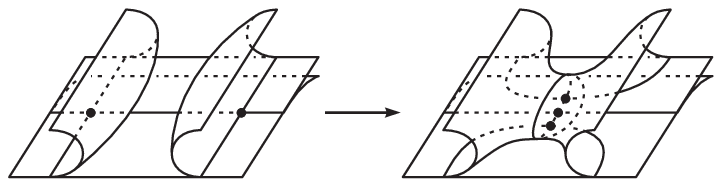}}
  \caption{$2\to 3$ bumping move.}\label{fig:bump}
\end{figure}

\noindent If $(P,b)$ is a branched spine and 
$P\to P'$ is either a  $2\to 3$ or a $0\to 2$ move, then it can be 
completed
(sometimes in a unique way, sometimes in two ways) to a {\it branched 
transit} $(P,b)\to (P',b')$. On the contrary, 
it could happen that a $3\to 2$ or $2\to 0$ inverse move 
is not ``branchable'' at all (see \cite[Ch.3]{20}). 
However we shall only use the ``positive'' moves. 
Note that otherwise it could stop any attempt to prove
the invariance of the quantum hyperbolic state sums, via any 
argument of ``move-invariance''. In Fig. \ref{fig:slide} and Fig. 
\ref{fig:bump} we show the $2\to 3$ branched transits. 

\begin{figure}
  \centerline{\psfig{figure=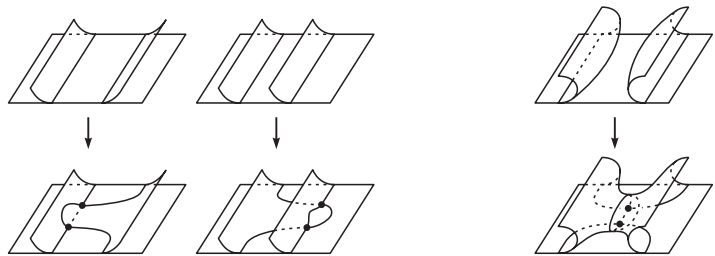}}
  \caption{branched lune-moves.}\label{fig:lune}
\end{figure}

\noindent In fact the list is not complete, but one can see in the figures 
all the essentially different behaviours, and easily complete the list
by applying evident symmetries. Following \cite{20} we can
distinguish two quite different kinds of branched transits: 
the {\it sliding moves} which actually preserve the combing, and the 
{\it bumping moves} which eventually change it. For the $0\to 2$ moves 
there is a similar behaviour (see Fig. \ref{fig:lune}).

\noindent However, we shall not
exploit this difference in the present paper. Note that the middle 
sliding move in Fig. \ref{fig:slide} corresponds dually to the triangulation 
move shown in Fig. \ref{w(e)}. 
   
\medskip

\noindent Given a distinguished triangulation $(T,H)$ of
$(W,L)$, the {\it first component} of $\mathcal{D}$ is just a 
branching $b$ of 
$T$.
If $(T,H)\to (T',H')$ is either a $2\to 3$, $0\to 2$ or $1\to 4$ move,
we have again some {\it branched transit}  $(T,H,b)\to (T',H',b')$; 
we have already described it in the first two cases. In the last case there are several ways to
take $b'$ which agree with $b$ on the edges ``already'' present in $T$.
Anyone of these ways is a possible transit.

\noindent By combining the existence of distinguished triangulations
with the above proposition we have:

\begin{lem} There exist almost-regular branched distinguished 
triangulations $(T,H,b)$ of $(W,L)$.
\end{lem}

\begin{remark} \label{rem2}{\rm {\bf Branching's r\^{o}le.} To understand the 
r\^{o}le of the branching, let us
go back to the ``classical'' hyperbolic ideal triangulations 
(see Remark \ref{rem1}). 
A hyperbolic ideal tetrahedron $\Delta$ with {\it ordered} 
vertices $v_0,\dots,v_3$ on the Riemann 
sphere $\mc \mathbb{P}^1 = \mc \cup \{ \infty \}$ 
is determined up to congruence by the {\it cross ratio} 
$z = [v_0:v_1:v_2:v_3]= \frac{v_2 - v_1}{v_2 - v_0} \frac{v_3 - 
v_0}{v_3 - v_1}$.
If the tetrahedron is positively oriented $z$ belongs to the upper half 
plane of $\mc$. Changing the ordering by {\it even} permutations produces
the cyclically ordered triple $(z,1/(1-z),(z-1)/z)$, which is the actual 
{\it modular} triple of $\Delta$. Changing the orientation corresponds 
to considering the complex conjugate triple.

\noindent A branching $b$ on the ideal triangulation allows one to choose,
in a somewhat {\it globally coherent} way, for each $\Delta$, a complex 
valued $z_{\Delta}$ in the corresponding triple. In the usual 
decoration of the edges of $\Delta$ by means of the elements
of the triple, $z_{\Delta}$ is associated to the first edge emanating 
from the base
vertex $v_0(\Delta)$, and so on, respecting the orders of the triple and of
the edges. 

\noindent Starting with a topological ideal triangulation of a
3-manifold $Z$, one usually tries to make it ``geometric'' 
(proving by the
way that $Z$ is hyperbolic) by solving some suitable system of 
{\it hyperbolicity} equations in $r_3$ complex
indeterminates, $r_3$ being the number of tetrahedra. Then, a branching 
allows us to specify \emph{which} system of equations, and this 
last is governed by a certain global coherence. 

\noindent Such a kind of global coherence
shall be important to control the behaviour of our state sums
up to decorated moves. On the other hand, it is not too surprising, in the 
spirit of the classical situation, that the value of the state sums 
shall not eventually depend on the branching of the decoration.}
\end{remark}

\subsection{Full cocycles}\label{full}

\noindent Recall that we are actually considering a triple $(W,L,\rho)$
where $\rho$ is an equivalence class of {\it flat} principal 
$B(2,\mc)$-bundles on $W$. Set $B= B(2,\mc)$.

\medskip

\noindent Let $(T,H,b)$ be a branched distinguished triangulation
of $(W,L)$. Let $[z]$ denote the equivalence class of a cellular 1-cocycle $z$ on $T$ via the usual 
equivalence relation {\it up to coboundaries}: since a (cellular) 0-cochain
$\lambda$ is a $B$-valued function defined on the vertices $V(T)$ of $T$, $z$ and $z'$ are
equivalent if they differ by the coboundary of some 0-cochain
$\lambda$. This means that for any $b$-oriented edge $e$ with ordered
end-points $v_0,v_1$, one has $z'(e)=\lambda (v_0)^{-1}z(e)\lambda
(v_1)$. We denote the quotient set by $H^1(T;B)$; recall that it can
be identified with the set of equivalence classes of flat
$B$-principal bundles on $W$.

\noindent The {\it second component} $z$ of 
$\mathcal{D}$ is a $B$-valued {\it full} 1-cocycle on $T$ representing 
$\rho$;  we write $[z] = \rho$. 
For each $b$-oriented
edge $e$ of $T$, $z(e)$ is an upper triangular matrix; we denote by 
$x(e)\in \mc$ the upper-diagonal entry of this matrix. ``Full'' means
that for each $e$, $x(e)\neq 0$.

\begin{remark} {\rm The {\it fullness} property strictly concerns the
cocycle $z$ and not its class $[z]$. It only depends on $T$ whether a 
given class
could be represented by full cocycles. Moreover, fullness does not
depend on the branching $b$. In fact we can define this notion by
using {\it any} arbitrary system of edge-orientations.}
\end{remark}

\noindent We refer to \S \ref{specialize} for more details and examples of
$B$-bundles. Here we simply recall that there are two distinguished 
abelian subgroups of $B$, and this fact induces some distinguished 
kinds of cocycles. 
They are: 

(1) the {\it Cartan} subgroup $C=C(B)$ of diagonal matrices; it is 
isomorphic to the multiplicative group $\mc^*$. The evident 
isomorphism $C\to \mc^*$ maps $A=(a_{ij})\in C$ to 
$a_{11}$; 

(2) the {\it parabolic} subgroup $Par(B)$ of matrices with double 
eigenvalue $1$; it is
isomorphic to the additive group $\mc$. The evident isomorphism maps 
$A=(a_{ij})\in Par(B)$ to $x=a_{12}$. 

\smallskip

\noindent Denote by $G$ any such subgroup. We get a map 
$H^1(T;G)\to H^1(T;B)$, where $H^1(T;G)$ is endowed with the natural
abelian group structure. Note that $H^1(T;Par(B))=H^1(T;\mc)$ 
is isomorphic to the ordinary (singular or de Rham) 1-cohomology of $W$. 

\noindent A complex valued {\it injective} function $u$ defined
on the set of vertices of a regular distinguished triangulation
$(T,H)$ of $(W,L)$ (as in the original set-up of \cite{1}) should be
regarded as a $Par(B)$-valued $0$-cochain.  Its coboundary $z=\delta
u$ is a basic example of a full cocycle representing $0\in
H^1(T;\mc)$ (whence the trivial flat bundle in $H^1(T;B)$).
\medskip

\noindent{\bf Existence and transit of full cocycles.} This is a 
somewhat delicate matter. Let $T$ be a triangulation of $W$, 
with {\it any} 
edge-orientation system $g$. Let $z$
be a 1-cocycle on $(T,g)$. Consider any transit $(T,g)\to (T',g')$ as
before. In the $2\to 3$ or $0\to 2$ cases (and their inverses), there
is a unique $z'$ on $(T',g')$ {\it which agrees with $z$ on the common
edges}. This defines a transit $(T,g,z)\to (T',g',z')$. Clearly
$[z]=[z']$.  

\noindent For $1\to 4$ moves $(T,g)\to (T',g')$, there is an infinite 
set of 
possible transits such that $[z]=[z']$ and the cocycles agree on
the common edges. Moreover, given one $(T,g,z)\to (T',g',z')$ $1 \to 4$ 
transit with $z$ a full cocycle, we can always turn $z'$ into an 
equivalent full $z''$, which differs from $z'$ by the coboundary of 
some $0$-cochain with support consisting of the new vertex $v$ of $T'$. 
Hence for $1\to 4$ moves there always exists an infinite set of {\it full} 
transits.

\noindent  Assume now that $z$ is full, and consider $2\to 3$ or $0\to 2$ 
moves (or their inverses). The trouble is that, in general, $z'$ is no
longer full. If $(T,H,b)\to (T',H',b')$ is a branched distinguished
transit and $z$ is full on $(T,b)$, then we have a completion to a
{\it full} transit $(T,H,b,z)\to (T',H',b',z')$ {\it only if} also the
final $z'$ is full. Otherwise we must stop.

\noindent However almost-regular triangulations have 
{\it generically} a good behaviour with respect to the existence and 
transit of full cocycles.

\begin{prop} \label{generique} Let $(T,g)$ be an almost-regular 
triangulation of $W$
endowed with any edge-orientations system $g$.
Let $(T,g)=(T_1,g_1)\to \dots \to (T_s,g_s)=(T',g')$ be any finite sequence
of $2\leftrightarrow 3$ or $0\leftrightarrow 2$ orientation-transits. 
Then there exists a
sequence $U_i$ of open dense sets of full cocycles on $T_i$
(in the natural subspace topology of $B^{r_1(T_i)}$), such that
$U_{i+1}$ is contained in the image of $U_i$ via the elementary transit
$T_i\to T_{i+1}$, and each class $\alpha \in H^1(T_i;B)$ can be represented
by cocycles in $U_i$. 
\end{prop}

\noindent As an immediate corollary (use the case $s=1$) we make a 
further step towards a solution of the {\it existence} problem:

\begin{cor} \label{fullOK}
Given a triple $(W,L,\rho)$, then a partially decorated (even almost-regular) distinguished triangulation $(T,H,b,z)$ do exists.
\end{cor}

\noindent {\it Proof of the proposition.} As $T$ is almost-regular, there
exists a regular triangulation $T''$ and a sequence 
$(T'',g'')\to \dots \to (T,g)$ of $2\leftrightarrow 3$ or 
$0\leftrightarrow 2$ orientation-transits. It is enough to prove the 
proposition for $T''$, so let us assume that $T$ 
is regular. The conclusion of the theorem holds for $s=1$. In fact we can 
suitably perturb any $z$ by the coboundaries of $0$-cochains, see 
the discussion above.
Now we remark that each elementary 
cocycle-transit $T_i\to T_{i+1}$ can be regarded as an algebraic bijective
map from the
space of 1-cocycles on $T_i$ to the space of 1-cocycles on $T_{i+1}$. 
The set of full 
cocycles for which the {\it full} elementary transit fails are contained
in a proper algebraic subvariety. So the conclusion follows, working by
induction on $s$.

\medskip

\noindent {\bf Full cocycles r\^{o}le.} Let $z$ be the full cocycle
of a decoration $\mathcal{D}$. For each edge $e$ denote by $t(e)$ the 
$(11)$-entry of $z(e)$; $x(e)$ is as before. Fix a determination $d$ of the
$N$'th-root holding for all entries of $z(e)$, for all $e$. 
Consider the Weil algebra $\mathcal{W}_N$, which is described in details in \S \ref{App}.
Then we can associate to each $e$ the $N$-dimensional standard
representation $r(e)$ of $\mathcal{W}_N$ characterized by the pair 
$(t(e)^{1/N},x(e)^{1/N})=(s(e),y(e))$.

\noindent This system of representations $\{r(e)\}$ satisfies the
following properties:

\smallskip

(1) Consider any $b$-oriented ``abstract'' 2-face $f$ of $T$. 
Starting from the base vertex $v_0(f)$, according to the orientation,
we find a cycle of edges $e_1,e_2,e_3$. The first two are positively
$b$-oriented. The last one has the negative orientation. Then 
the cocycle condition on $z$ implies that 

\smallskip

\emph{$r(e_3)$ is, up to isomorphism, the unique irreducible summand of the 
representation $r(e_1)\otimes r(e_2)$.} 

\smallskip

(2) Consider, in particular, the 2-face $f$ opposite to
the base vertex $v_0(\Delta)$, in any tetrahedron $\Delta$ of $T$.
For each $e_i$ of $f$, let us denote by $e_{o(i)}$ the opposite edge in 
$\Delta$.
Then the $y$-components of the representations $r(e_i),r(e_{o(i)})$ satisfy 
the following {\it Fermat} relation:
$$ y(e_3)^Ny(e_{o(3)})^N = y(e_1)^Ny(e_{o(1)})^N
   + y(e_2)^Ny(e_{o(2)})^N \ .$$

(3) The same conclusions hold for any other branching
$b'$ on $T$, as they only depend on the cocycle condition.
\medskip

\noindent After a full transit (once also the charge-transit shall
be ruled out), we observe that: 

\smallskip

\noindent {\it The system of representations $\{r(e)\}$ varies 
exactly in the way one needs in order to apply the fundamental
algebraic identities (the pentagon, orthogonality, and bubble relations of the Appendix) satisfied 
by the quantum-dilogarithm c-6j-symbols.}

\smallskip

\noindent One could ask if systems of representations $r=\{r(e)\}$ 
verifying the above properties are more general than the one 
obtained starting from full cocycles. In fact, setting $t(e)=s(e)^N$ and
$x(e)=y(e)^N$ one obtains a full cocycle $z$. In a sense the system
$\{r(e)\}$ can be considered as a sort of {\it reduction} mod $N$ of $z$. 
In order to  study any ``classical'' limit,
when $N\to \infty$, of our state sums invariants, it seems quite 
appropriate to consider the
reductions of the {\it same} cocycle.   

\medskip  

\subsection{ Charges} \label{charges}

\noindent The ``classical'' source of the charges $c$ in the decorations $\mathcal{D}$ clearly emerges from Neumann's work on the Cheeger-Chern-Simons classes of hyperbolic $3$-manifolds and scissors congruences of hyperbolic polytopes. Since there is a wide literature on this subject (see the references in \cite{13}), we shall only report a few details. 

\medskip

\noindent{\bf Refined Scissors Congruence.}
From the work of \cite{23} and \cite{26}, we know that the 
volume of any oriented hyperbolic
3-manifold $\mathcal{F}$ has a deep analytic relationship with
another geometric invariant, the Chern-Simons invariant $CS(\mathcal{F})$. Recall that $CS(.)$ is a $\mr/2\pi^2\mz$-valued invariant defined for any oriented compact Riemannian 3-manifold, and that its definition can also be extended, with value in $\mr/\pi^2\mz$, to \emph{non-compact} complete and finite volume hyperbolic 3-manifolds \cite{27}. Consider then
$$VCS(\mathcal{F}) = Vol(\mathcal{F})+iCS(\mathcal{F}).$$ 

\smallskip

\noindent It turns out that $VCS(\mathcal{F})$ (for 
any finite volume $\mathcal{F}$), actually depends on a weaker
sub-structure of the full hyperbolic structure, called the {\it
refined scissors congruence class} and denoted by 
$\widehat{\beta} (\mathcal{F})$ \cite{13,14}; it is orientation 
sensitive and takes values
in the {\it extended Bloch group} $\widehat{\mathcal{B}} (\mc)$. 
If $\mathcal{F}$ is non-compact, this class may be represented with the 
help of
the usual geometric ideal triangulations $T$; if $\mathcal{F}$ is
compact, we assume here that it is obtained by {\it
hyperbolic} Dehn surgery on some non-compact one, so that we have some
deformed geometric ideal triangulations $T$ of 
$\mathcal{F}\setminus \mathcal{G}$, where ${\mathcal{G}}=\{L_i\}$ is a 
finite set 
of simple closed geodesics in $\mathcal{F}$. Then, by abuse of 
notations, in both cases we shall
denote these special decompositions of $\mathcal{F}$ by $T$ and call them 
``ideal triangulations of $\mathcal{F}$''. 

\smallskip

\noindent 
For the sake of clarity, let us present the construction of
$\widehat{\beta}(\mathcal{F})$. It relies heavily on the functional
properties of the classical {\it Rogers dilogarithm}, defined on $\mc
\setminus \left( ]-\infty,0[ \cup ]1,\infty[ \right)$ by:
$$ \mathcal{R}(z) = \frac{1}{2} \log(z)\log(1-z) - \int_{0}^{z}  
\frac{\log (1-t)}{t} dt  - \pi^2/6 \ ,$$
\noindent and in particular on its five-term identity, which reads:
\begin{eqnarray}\label{five-term}
\mathcal{R}(x) + \mathcal{R}(y) - \mathcal{R}(xy) 
= \mathcal{R}\left( \frac{x(1-y)}{1-xy}\right) 
+ \mathcal{R} \left( \frac{y(1-x)}{1-xy} \right)\ .
\end{eqnarray}
\noindent One then starts with a four components non connected covering of 
$\mc \setminus \{0,1\}$: 
$$\widetilde{\mc} =X_{00}\cup X_{10}\cup X_{01}\cup X_{11} \ .$$

\noindent It can be regarded as the Riemann surface for the collection 
of all branches of the functions $(\log(z)+p\pi i, -\log(1-z)+q\pi i)$
on $\mc\setminus \{0,1\}$, $(p,q)\in \mz \times \mz$. If $P$ is
obtained by splitting $\mc\setminus \{0,1\}$ along the rays
$(-\infty,0)$ and $(1,\infty)$, $\widetilde{\mc}$ is a suitable
identification space from $P\times \mz \times \mz$: namely,
$X_{\epsilon_1,\epsilon_2} = \{ [z,p,q] \in \widetilde{\mc} \ \vert \
p \equiv \epsilon_1 {\pmod 2},q \equiv \epsilon_2 {\pmod 2}\}$. The
map
$$l([z,p,q])=(\log (z) + p\pi i, 
-\log (1-z)+q\pi i, \log (1-z)-\log (z) -
(p+q)\pi i)$$

\noindent 
is well-defined on $\widetilde{\mc}$, and it gives an identification
between $\widetilde{\mc}$ and the set of triples of the form
$$(w_0,w_1,w_2)=(\log (z) + p\pi i, \log (z') +q\pi i, \log (z'') 
+ r\pi i)$$
\noindent with
$$ p,q,r\in \mz\ \ {\rm and}\ \ \sum_i w_i = 0 \ ,$$
\noindent and for some determination of the $\log$ function, where 
$(z,z',z")$ is a modular triple for an ideal tetrahedron $\Delta$; in
other words we adjust its dihedral angles by means of multiples of
$\pi$ so that the resulting angle sum is zero. Such a triple is called a {\it
combinatorial flattening} of the ideal tetrahedron. So
$\widetilde{\mc}$ can be regarded as the set of these combinatorial
flattenings, for all ideal tetrahedra of the hyperbolic space
$\mathbb{H}^3$.

\noindent The Rogers dilogarithm lifts to an analytic function:
$$R: \widetilde{\mc} \to \mc/\pi^2\mz$$ 
$$R([z,p,q])=\mathcal{R}(z)+\frac{i\pi}{2}(p\log (1-z)+q\log (z))\ ,$$
\noindent and it extends to the free 
$\mathbb{Z}$-module $\mathbb{Z}[\widetilde{\mc}]$.  Moreover, one may
lift the classical five-term identity satisfied by the Rogers
dilogarithm to the function $R$, provided that we pass to some
quotient of $\mathbb{Z}[\widetilde{\mc}]$. There is such a natural
maximal one $\widehat{\mathcal{P}} (\mc)$, which is called the
\emph{extended pre-Bloch group}. In this way, we can turn $R:
\widehat{\mathcal{P}} (\mc) \to \mc/\pi^2\mz$ into a homomorphism.

\noindent  Take our
finite volume oriented hyperbolic 3-manifold $\mathcal{F}$; for any
ideal triangulation $T$ of $\mathcal{F}$, one can define
$$\widehat{\beta} (T) = \sum_{\Delta \in T} \eta_{\Delta} 
[z,p,q](\Delta) \in \widehat{\mathcal{P}} (\mc),$$
\noindent where $\eta_{\Delta}\in \{1,-1\}$ depends on the orientation
of $\Delta \in T$. When $\mathcal{F}$ is non-compact 
$\widehat{\beta} (T)$ just represents the (refined) scissors congruence
class  $\widehat{\beta} (\mathcal{F})$ of  $\mathcal{F}$. 
In the compact case, how to explicitly represent 
$\widehat{\beta} (\mathcal{F})$ is a subtler fact; anyway the above  
$\widehat{\beta} (T) =  \widehat{\beta} (\mathcal{F},\mathcal{G})$ 
is adequate to represent
a ``classical'' counterpart of our state sums, which in general depend
also on the link and not only on the ambient manifold.

\begin{remark} \label{inop} {\rm Consider the map:
$$\widehat{\delta} :\widehat{\mathcal{P}} (\mc)\to \mc
\wedge_{\mathbb{Z}}\mc$$
$$[z,p,q] \mapsto (\log(z) + ip\pi) \wedge (-\log(1-z) + iq\pi).$$
\noindent It lifts the Dehn invariant of the classical scissors
congruence \cite{13}. Then $\widehat{\mathcal{B}} (\mc)=
Ker\widehat{\delta}$ is called the {\it extended Bloch group}, and one
can show that $\widehat{\beta} (\mathcal{F})$ or 
$\widehat{\beta} (\mathcal{F},\mathcal{G})$ are in
$\widehat{\mathcal{B}} (\mc)$.}
\end{remark} 

\noindent The classes $\widehat{\beta} (\mathcal{F})$ or 
$\widehat{\beta}(\mathcal{F},\mathcal{G})$ may be seen 
as representatives of the fundamental class of $\cal F$ in
the (discrete) homology group (\cite{14,15})
$$H_3^{\delta}(PSL(2,\mc);\mz),$$
\noindent which itself maps surjectively onto
$\widehat{\mathcal{B}} (\mc)$ \cite{14} (see the definition in Remark
\ref{inop}). The relations in $\widehat{\mathcal{P}} (\mc)$
express the $2 \to 3$ move for ideal hyperbolic tetrahedra endowed
with combinatorial flattenings, and then they give the independence of
$\widehat{\beta} (T)$ from $T$.

\bigskip

\noindent The relations in $\widehat{\mathcal{P}} (\mc)$ allow to have a 
global
control on the values of the combinatorial flattenings, and this is
made easier by the presence of a branching $b$ of $T$. For instance,
note that a, let us say, $X_{00}$-flattening of a $\Delta$ presupposes
a choice of opposite edges; changing this choice turns it into a
$X_{10}-$ or $X_{01}-$flattening.

\noindent Anyway, the point here is that part of this global coherence is 
expressed in terms of relations which must be satisfied by the 
{\it integral} components $(p,q)$ of the flattenings $\{[z,p,q](\Delta)\}$, 
which then give, by definition, an {\it integral charge} on $T$. 
These relations are strongly reminiscent of ``the''
system of {\it hyperbolicity} equations (for instance the one specified by 
the branching) which is satisfied by the {\it arguments} of 
the modular triples of the ideal triangulation. The charge $c$ of 
a decoration $\mathcal{D}$ shall be, essentially, the {\it reduction} 
mod $N$ of an integral charge. In \S \ref{GR} we shall return on
these scissors congruence classes and on their relationship
with the $VCS$-invariant.   

\medskip

\noindent{\bf Integral charges.} We shall now give the formal 
definition of the
{\it integral charges} in our own setting. It is a straightforward 
adaptation
of the integral charges of the previous paragraph. 

\noindent Let $(T,H)$ be a distiguished triangulation of $(W,L)$.
With the notations of Definitions \ref{adapte} and \ref{type}, let us
assume for simplicity that the associated {\it simple} spine
$\widetilde{P}$ of $M$ is in fact a standard one, as in the proof of
Lemma \ref{existadapte}. The following discussion could be adapted
also to the case when $\widetilde{P}$ is merely simple; anyway, we
could also add the $\widetilde{P}$ standardness assumption to our
set-up without any substantial modification in all our arguments.

\smallskip

\noindent We know that the truncated tetrahedra of $T(\widetilde{P})$ 
induce a 
triangulation
$\tau$ of $\partial M$. 
Let $s$ be an oriented simple closed curve on $\partial M$ in 
general
position with respect to $\tau$. We say that $s$
{\it has no back-tracking} with respect to $\tau$ if it never departs a 
triangle of $\tau$ across the same edge by which it entered. Thus each time 
$s$ passes through a triangle, it selects the vertex between the 
entering and departing edges, and gives it a sign $\in\{1,-1\}$, 
according as it goes past this vertex positively
or negatively with respect to the boundary orientation.

\noindent Let $s$ be a simple closed curve in $M$ in general position
with respect to the ideal triangulation $T(\widetilde{P})$. We say
that $s$ {\it has no back-tracking} with respect to $T(\widetilde{P})$
if it never departs a tetrahedron of $T(\widetilde{P})$ across the
same 2-face by which it entered.  Thus each time $s$ passes through a
tetrahedron, it selects the edge between the entering and departing
faces.

\noindent Fix in each component of $\partial M$ one of the 
meridians in the boundary of $Q$, $m_i$ say, and a simple closed curve
$l_i$ in general position with respect to $\tau$, which intersects
$m_i$ transversely at one point. Orient $m_i$ and $l_i$. We can assume
that these curves have no back-tracking with respect to $\tau$. It is
clear that any class in $H_1(\partial M; \mz)$ may be represented in
the set $\{m_i,l_i\}$ of isotopy classes of curves without
back-tracking with respect to $\tau$ and generated by $m_i$ and
$l_i$. In the same manner, any class in $H_1( M; \mz/2\mz)$ may be
represented up to isotopy by a curve in $M$ without back-tracking with
respect to $T(\tilde{P})$.

\noindent Denote by $E_{\Delta}(T)$ the set of all edges of all ``abstract'' 
tetrahedra of $T$;
there is a natural map $\epsilon :E_{\Delta}(T) \to E(T)$.

\begin{defi} \label{defcharges}  {\rm An {\it integral charge} 
on $(T,H)$ is a map
$$c': E_{\Delta}(T) \to \mz$$ which satisfies the following
properties:

\smallskip

\noindent (1) For each 2-face $f$ of any abstract $\Delta$ with edges 
$e_1,e_2,e_3$,   
$$\sum_i c'(e_i) = 1 \ ,$$
\noindent for each $e\in E(T)\setminus E(H)$,
$$ \sum_{e'\in \epsilon ^{-1}(e)}c'(e') = 2\ ,$$
\noindent for each $e\in E(H)$,
$$ \sum_{e'\in \epsilon ^{-1}(e)}c'(e') = 0\ .$$

\noindent (2) Let $s$ be a curve on $\partial M$ which has no
back-tracking with respect to $\tau$. Each time $s$ enters a triangle
of $\tau$, $c'$ associates in a natural way an integer to the selected
vertex; multiply this integer by the sign of the vertex and take the
sum $\sigma(s)$ of these signed integers. Then for every $s\in
\{m_i,l_i\}$,
$$\sigma(s)=0\ .$$

\noindent (3) Let $s$ be any curve which has no back-tracking with 
respect to $T(\widetilde{P})$. Each time $s$ enters a tetrahedron of 
$T(\widetilde{P})$,
$c'$ associates in a natural way an integer to the selected edge. Take
the sum $\alpha(s)$ of these integers. Then, for each $s$,
$$\alpha(s) \equiv 0\ \ {\rm mod}\ 2\ .$$}  
\end{defi}

\begin{defi} {\rm The {\it third} component $c$ of any 
decoration $\mathcal{D}$,
called the {\it charge}, is of the form:
$$c \equiv \frac{c'}{2}\ {\rm mod} \ N,$$
where $1/2 = p +1 \in \mz/N\mz$ (with the notations of the Appendix), and
$c'$ is an integral charge on $(T,H)$.}
\end{defi}

\begin{remark}{\rm By conditions (1), $c$ satisfies, in particular,
the charge requirements of the quantum data, see Proposition
\ref{symmetry}. By the way, note the importance of $N$ being odd in
the present definition.  Conditions (2) and (3) are in fact purely
homological; in particular, as already said, their validity does not
depend on the particular choice of $m_i$ and $l_i$ on each component
of $\partial M$.}
\end{remark}

\noindent{\bf Charge's existence and transit.} The existence of integral
charges is obtained by just rephrasing the proof of Theorem 2.4.(i)
(that is of Lemma 6.1) in \cite{28}. The only modification is in
considering $T$ as an ideal triangulation of $W \setminus rD^3$ (see Remark \ref{trig}); if we set the sum of the charges around 
the edges of $H$
to be equal to $0$, the existence follows via the same arguments.

\noindent Also Theorem 2.4.(ii) \cite{28} shall be important in order to prove 
the state sum charge-invariance. Let us first describe qualitatively
this result. Let $r_0$ and $r_1$ be respectively the number of vertices and edges of $T$; an easy computation with the Euler characteristic shows that
there are exactly $r_1 - r_0$ tetrahedra in $T$ (see Proposition
\ref{invproj}). The first condition in Definition \ref{defcharges} (1)
says that there are only two independent charges on the edges of each
abstract $\Delta$, and given a branching $b$ on $T$ there is a
preferred such ordered pair $(c_1^{\Delta},c_2^{\Delta})$. Set
$c_1^{\Delta}=w_1^{\Delta},c_2^{\Delta}=-w_2^{\Delta}$ (then
$c_3^{\Delta}=-w_1^{\Delta}+w_2^{\Delta}+1$). Then, there is a
canonical way to write down an integral charge on $(T,b)$ as a vector
in $\mz^{2(r_1-r_0)}$ with first components $w_1^{\Delta},\ \Delta \in T$,
and then $w_2^{\Delta},\ \Delta \in T$.
\noindent Theorem 2.4.
(ii) in \cite{28} says that integral charges lie in a specific
sublattice of $\mz^{2(r_1-r_0)}$:

\begin{prop} \label{lattice} 
There exist determined $w(e) \in \mz^{2(r_1-r_0)}$, $e \in T$, such that, given any
integral charge $c'$, all the other integral charges $c''$ are of the
form
$$c''=c'+\sum_e \lambda_ew(e) $$
\noindent 
where the second addendum is an arbitrary integral linear combination
of the $w(e)$.
\end{prop}

\begin{figure}[ht]

\begin{center}

\scalebox{1}{\input{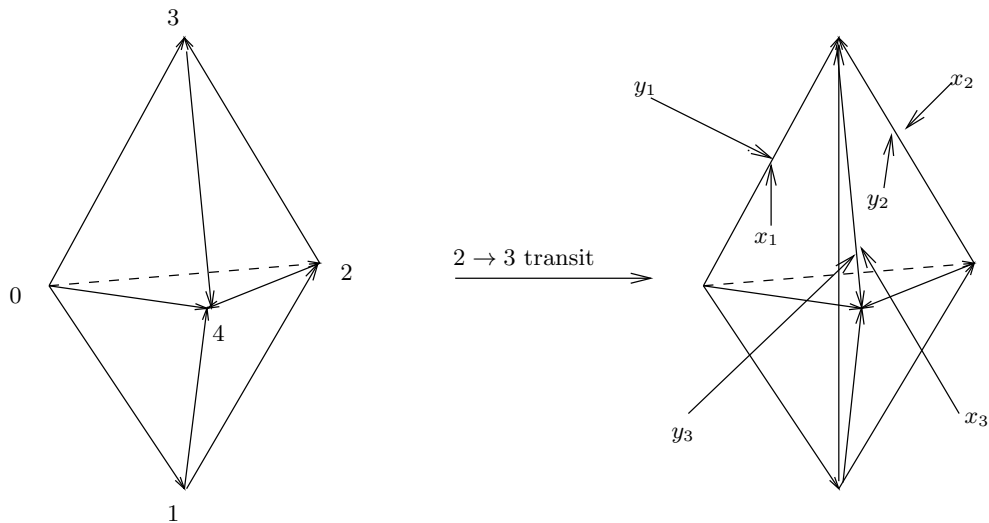}}

\end{center}

\caption{\label{w(e)} $2 \to 3$ charge transits 
are generated by Neumann's vector $w(e)$}

\end{figure}

\noindent 
The vectors $w(e)$ have the following form. For each tetrahedron
$\Delta \in (T,b)$ glued along a specific $e$, define $r_1^{\Delta}(e)$ and
$r_2^{\Delta}(e)$ as the coefficients in $c'(\epsilon^{-1}(e))$, when
written in terms of $(w_1^{\Delta},w_2^{\Delta})$. Then $w(e)$ is the
vector in $\mz^{2(r_1-r_0)}$ with first components $r_2^{\Delta},\ \Delta
\in T$, and then $-r_1^{\Delta},\ \Delta \in T$. For instance, in the
situation described on the right of Fig. \ref{w(e)} (where the ordering
of the tetrahedra is induced by the ordering of the vertices), we
easily see that $w(e) = (-1,1,-1,1,0,1)^t$.

\medskip

\medskip

\noindent Next, we describe the {\it transit of integral charges}. 
The transit of $\mathcal{D}$-charges shall be obtained 
by reduction mod $N$. Let $(T,H)\to (T',H')$ be a $2\to 3$ move. Let
$c$ be an integral charge on $(T,H)$, and $e$ the edge that
appears. Consider the two ``abstract'' tetrahedra $\Delta_1,\Delta_2$
of $T$ involved in the move. They determine a subset $\widetilde{E}
(T)$ of $E_{\Delta}(T)$.  Denote by $\widetilde{c}$ the restriction of
$c$ to $\widetilde{E} (T)$.  Let $c'$ be an integral charge on
$(T',H')$. Consider the three ``abstract'' tetrahedra of $T'$ involved
in the move. So we have $\widetilde{E} (T')$ and $\widetilde{c}'$ with
the clear meaning of the symbols. Denote by $\widehat{E} (T)$ the
complement of $\widetilde{E} (T)$ in $E_{\Delta}(T)$ and $\widehat{c}$
the restriction of $c$.  Do similarly for $\widehat{E} (T')$ and
$\widehat{c}'$. Clearly $\widehat{E} (T)$ and $\widehat{E} (T')$ can
be naturally confused. The following lemma is the key point of the
charge-invariance.

\begin{lem-defi} \label{charge-transit} 
{\rm We have a {\it charge-transit} $(T,H,c)\to (T',H',c')$ if:

\smallskip

(1) For each ``common edge'' $e_0 \in 
\epsilon_T(\widetilde{E} (T))\cap \epsilon_{T'}(\widetilde{E} (T'))$ of $T$ and $T'$,
$$ \sum_{e'\in \epsilon_T^{-1}(e_0)}\widetilde{c} (e')=
 \sum_{e"\in \epsilon_{T'}^{-1}(e_0)}\widetilde{c}' (e")\ .$$

(2) $\widehat{c}$ and $\widehat{c}'$ agree on $\widehat{E} (T)=
    \widehat{E} (T')$.

\noindent Moreover, if $c$ 
is an integral charge on $(T,H)$ and a function $c'$ satisfies these
two conditions, then $c'$ is an integral charge.}
\end{lem-defi}

\noindent 
Charge-transits for $0\to 2$ and $1\to 4$ moves are defined in a
similar way; in particular, notice that charges are defined without
any reference to branchings, and that non-branched $0 \to 2$
moves may be obtained from $2 \to 3$ and $3 \to 2$ moves. Moreover, for $1 \to 4$ moves, one must add the condition $(3)$ that the sum of the charges is $0$ around the chosen two new edges of $H$, and $2$ around the others that appear.

\medskip

\noindent 
\emph{Proof.} For this definition to make sense, we have to show that such
transits actually define integral charges, and by (2) we can restrict
our attention to Star($e,T$). Now, assuming (1), the identity $\sum_i
c(e_i) = 1$ on $f = \Delta_1 \cap \Delta_2$ in Definition
\ref{defcharges} (1) is equivalent to $ \sum_{e'\in \epsilon
^{-1}(e)}c'(e') = 2$ (compare with formula (\ref{e0}) in \S \ref{App}).  Hence
\ref{defcharges} (1) for $c'$ is clearly verified. One can also easily show that $c'$ satisfy \ref{defcharges} (2)-(3), for they are homological conditions. This proves our claim.

\noindent 
Consider now the situation of Fig. \ref{w(e)}; recall that the
branching is only used to write down explicitly the set of charge
coefficients.  It is easy to verify that, assuming (1) and (2), the
$\mz$-vector space of solutions $(x_1,y_1,\ldots, x_3,y_3)$ to the
system of linear equations defining $c'$ from $c$ is generated by
$w(e) = (-1,1,-1,1,0,1)^t$. The integral charges on $T'$ may then only
differ by a $\mz$-multiple of $w(e)$, thus concluding the proof.
 
\subsection{Total decorations} \label{totdeco} 

\noindent Summing up: 

\medskip

(1) we dispose now of a complete
solution of the {\it existence} problem.

\begin{teo} For any $(W,L,\rho)$ there exist (even almost-regular) 
decorated distinguished triangulations 
$\mathcal{T}=(T,H,\mathcal{D}=(b,z,c))$.
\end{teo}

(2) We have clearly specified what is a $2\to 3$, $0\to 2$
or $1\to 4$ {\it transit} $\mathcal{T}\to \mathcal{T}'$ between
decorated distinguished triangulations; sometimes we refer to
them  as {\it decorated} moves.
\medskip

(3) We have pointed out how the decorations are reminiscent of the 
{\it refined scissors congruence} classes of hyperbolic 3-manifolds,
represented by means of geometric ideal triangulations.

\section{The invariant $K_N(W,L,\rho)$}\label{WL}

This section contains the construction of the state sums, the proof of
their invariance and  the proof of some of their basic properties.

\subsection{State sum transit-invariance}\label{state-sum1}

\noindent
 Let $\mathcal{T}=(T,H,\mathcal{D}=(b,z,c))$ be as above. The state
sum is defined by:
\begin{eqnarray} \label{invariant}
K_N(\mathcal{T},d)=\left( N^{-r_0}\sum_{\alpha} \prod_{\Delta}
t^{\Delta}(\mathcal{D},\alpha,d)\prod_{e\in E(T)\setminus E(H)}
x(e)^{-2p/N} \right)^N\ .
\end{eqnarray}
\noindent In fact, denoting by $F(T)$ the set of 2-faces of $T$, each 
{\it state} 
$$\alpha:F(T)\to \mz/N\mz$$ 
\noindent together
with $\mathcal{D}$ induce a {\it complete} decoration on each
tetrahedron $\Delta$, and hence an associated quantum-dilogarithm
c-$6j$-symbol $t^{\Delta}(\mathcal{D},\alpha,d)$ (see \S \ref{App} for
details). Recall that $x(e)\neq 0$ is the upper-diagonal entry of
$z(e)$ and that $d$ is a determination of the $N$'th-root holding for
all entries of $z(e)$ (see \S \ref{full}).  Consider any decorated transit $\mathcal{T}\to
\mathcal{T}'$. Assume that $d$ holds for both $\mathcal{T}$,
$\mathcal{T}'$. Then we have:

\begin{prop} \label{transitinv} The state sum is transit invariant:
$$K_N(\mathcal{ T},d)=K_N(\mathcal{ T}',d)\ .$$
\end{prop}

\noindent {\it Proof.} In some sense this is the main achievement of \cite{1}.
The equality is a consequence of the fundamental algebraic identities
satisfied by the quantum-dilogarithm c-$6j$-symbols: the pentagon ($2\to
3$), othogonality ($0\to 2$), and bubble relations. The precise
statements are given respectively in Propositions \ref{EP}, \ref{orth}
and \ref{vertex}.

\subsection{ State sum total invariance}\label{state-sum2}

\noindent We can finally state our first main result.

\begin{teo} \label{teo1} Let $\mathcal{ T}=
(T,H,\mathcal{ D}=(b,z,c))$ be a decorated distinguished triangulation
of $(W,L)$ such that $[z]=\rho$.  Then, for any odd positive integer
$N$, $K_N(W,L,\rho)=K_N(\mathcal{ T},d)$ is a well-defined
invariant. In particular, the original Kashaev set-up corresponds to
the case of the trivial flat bundle $\rho$. These are purely
topological invariants $K_N(W,L)$.
\end{teo}

\noindent The following proposition is formally contained in the
statement of the theorem.

\begin{prop} \label{propo1} 
(a) Let $\mathcal{ T}=(T,H,\mathcal{ D})$ be a decorated distinguished
triangulation of $(W,L)$. Let $d$ and $d'$ be two determinations of
the $N$'th-root holding for $\mathcal{T}$. Then $K_N(\mathcal{T},d)=
K_N(\mathcal{ T},d')$.  That is, the root determination has
no uninfluence, and we can omit indicating it.

\noindent 
(b) Let $\mathcal{D}$ and $\mathcal{D}'$ be two decorations on
$(T,H)$, giving $\mathcal{ T}$ and $\mathcal{ T}'$ . Then
$K_N(\mathcal{ T})=K_N(\mathcal{ T}')$.
\end{prop}

\noindent 
{\it Proof of Theorem \ref{teo1} assuming Proposition \ref{propo1}}
Let us fix a model of $W$, with $L \subset W$ considered up to
isotopy. We have to prove that the state sum value does not depend on
the decoration $\mathcal{T}$ used to compute it. Consider $\mathcal{
T}$ and $\mathcal{ T}'$, and assume first that the corresponding $T$
and $T'$ are {\it almost-regular}.  Up to some $1\to 4$ moves (use
transit-invariance) we can assume that $T$ and $T'$ have the same
vertices and that the corresponding spines of $W$ adapted to $L$ have
the same type and coincide along $L$. Then, let us apply Proposition
\ref{chemintype} to $T$ and $T'$; we find $T''$; as Proposition
\ref{chemintype} uses only ``positive'' moves, the branchings
transit. Moreover, we may choose the full cocycles on $T$ and $T'$
{\it generically} (use \ref{generique} and \ref{transitinv} a first
time) so that one realizes a full cocycle-transit. Finally, the
charges also transit by Lemma \ref{charge-transit}.  Summing up, we
have two decorated transits $\mathcal{ T}\to \mathcal{ T}_1$,
$\mathcal{ T}'\to \mathcal{ T}_2$, where at the final step one
possibly has different decorations of the same $(T'',H'')$. The
transit-invariance gives us
$$K_N(\mathcal{ T})= K_N(\mathcal{ T}_1)$$ 
\noindent and 
$$K_N(\mathcal{ T'})= K_N(\mathcal{ T}_2)\ .$$
\noindent Finally  Proposition \ref{propo1} gives us 
$$K_N(\mathcal{ T}_1)=K_N(\mathcal{ T}_2)\ .$$
\noindent We are done in the almost-regular setting. To finish it is
enough to show that any decorated distinguished triangulation transits
to an almost-regular one. This is easily achieved by means of 
$1\to 4$ moves. Since one can identically pull back the decorations via (PL) 
homeomorphisms of the triples $(W,L,\rho)$, the proof is complete.

\bigskip

\noindent {\it Proof of Proposition \ref{propo1}}. 
Let us fix as above a model of $W$, with $L \subset W$ considered up
to isotopy.

\medskip

\noindent{\bf $N$'th-root determination invariance.} We have to prove (a)
of \ref{propo1}. This is a consequence of 
the fact that the functions $h$ and $\omega$
in the quantum 
dilogarithm c-$6j$-symbols are homogeneous of degree $0$; see \S \ref{App}.

\medskip
  
\noindent{\bf Branching-invariance.} 
Let $\mathcal{T}$ and $\mathcal{ T}'$ be two decorated distinguished
triangulations of $(W,L)$ such that they only differ by the branchings
$b$, $b'$. We have to prove that $K_N(\mathcal{ T})=K_N(\mathcal{
T}')$.  This is done in Lemma \ref{indbranch}. It describes
$K_N(\mathcal{ T})$ as the $N$'th-power of a \emph{weighted trace} of
an operator defined on some tensor product of standard representations
of $\mathcal{ W}_N$, determined by the full cocycle $z$ and the
determination $d$.

\medskip

\noindent{\bf Charge-invariance.} Let $\mathcal{T}$ and $\mathcal{ T}'$ be
decorated distinguished triangulations of $(W,L)$ such that the two
decorations only differ by the charges $c$, $c'$. We have to prove
that $K_N(\mathcal{ T})=K_N(\mathcal{ T}')$. We may assume that $T$ is
almost-regular. We shall use the transit-invariance: so we probably
have to generically change the full cocycles in order to guarantee
full transits. Anyway, this is not a trouble because, by continuity,
it is enough to prove the present statement for full cocycles
arbitrarily close to the one of $\mathcal{T}$. Let us fix an integral
charge on $(T,H)$ which we denote again $c$, of which the
$\mathcal{D}$-charge is the reduction mod $N$. Fix any edge $e$ of
$T$. Consider all the charges of the form (we use the notations of
Proposition \ref{lattice})
$$ c'= c+ \lambda w(e),\ \ \lambda \in \mz \ .$$
\noindent Let us denote this set of charges $C(e,c,T)$. It is the set
of charges which differ from $c$ only on Star$(e,T)$.  Thanks to
Proposition \ref{lattice}, it is enough to prove the charge-invariance
when $c'$ varies in $C(e,c,T)$. In this way we have somewhat
``localized'' the problem. The result is an evident consequence of the
following facts:

(1) Let $\mathcal{ T}\to \mathcal{ T}''$ be any $2\to 3$
 transit such that $e$
is a common edge of $T$ and $T''$. Then the result holds for $C(e,c,T)$
iff it holds for $C(e,c'',T'')$.

(2) There exists a sequence of $2\to 3$ transits 
$\mathcal{ T}\to \dots \to \mathcal{ T}''$ 
such that $e$ persists at each step,
and Star$(e,T'')$ is like the final configuration of a $2\to 3$ move
with $e$ playing the role of the central common edge of the 3
tetrahedra.

(3) If Star$(e,T)$ is like Star$(e,T'')$ as above, then the result holds for 
$C(e,c,T)$.

\noindent Property (1) is a 
consequence of the transit-invariance (Proposition \ref{transitinv}),
because $C(e,c,T)$ transits to $C(e,c'',T'')$ by the first claim of Lemma
\ref{charge-transit}.

\noindent To prove (2) it is perhaps easier to think, for a while, 
in dual terms. Consider the dual region $R=R(e)$ in $P(T)$. There is
a natural notion of {\it geometric multiplicity} $m(R,a)$ of
$R$ at each edge $a$ of $P$, and $m(R,a)\in \{0,1,2,3\}$.  We say that
$R$ is {\it embedded} in $P$ iff for each $a$, $m(R,a)\in \{0,1 \}$.
To describe the final configuration of $e$ in $T''$ is equivalent to
saying that the dual region is an embedded {\it triangle}. Each time $e$
has a {\it proper} (i.e. with two distinct vertices) edge $a$ with
$m(R,a)\in \{2,3 \}$, the $2\to 3$ move at $a$ ``replaces'' $a$ with
new edges $a'$ with $m(R,a')< m(R,a)$.  If $R$ has only loops with bad
multiplicity, a suitable $2\to 3$ move at a proper edge of $P(T)$ with
a common vertex with the loop ``replaces'' the loop with proper
edges. By induction we get that, up to $2\to 3$ moves, we can assume
that $R$ is an embedded polygon. Look now at the dual situation. We
possibly have more than 3 tetrahedra around $e$. It is not hard to
reduce the number to 3, via some further $2\to 3$ moves.

\noindent Property (3) is almost an immediate consequence of the $3\to 2$ transit.
In fact, given any $2\to 3$ charge-transit $(T,c)\to (T',c')$, we know
by the second claim of Lemma \ref{charge-transit} that all 
the other charges $c''$,
varying the transit $(T,c)\to (T',c'')$, exactly make
$C(e,c',T')$. But there is a little subtlety: in
general, the branching $b$ does not transit during a $3\to 2$
move. Anyway, we can modify the branching $b$ on the 3 tetrahedra
around $e$ in such a way that the $3\to 2$ move becomes
branchable. So we have on $T$ the original branching $b$ and another
system of edge-orientations $g$. Then apply Proposition
\ref{brancheable} to $(T,g)$; we find $T'$ with two branchings: $b'$
by the transit of $b$ (recall that \ref{brancheable} uses only
positive moves) and the branching $b''$ over $g$. Note that, thanks to
the actual proof of \ref{brancheable}, $e$ persits and Star$(e,T')=$
Star$(e,T)$. Moreover we have a charge-transit $(T,c)\to (T',c')$ with
$c$ and $c'$ which agree on Star$(e,T')$. So, using the
branching-invariance, we may assume that the $3\to 2$ move is
branchable, and the charge-invariance is thus proved.

\medskip 

\noindent{\bf Cocycle-invariance.} Let $\mathcal{T}$ and $\mathcal{ T}'$ be
two decorated distinguished triangulations of $(W,L)$ which only
differ by the full cocycles $z$, $z'$, $[z]=[z']$.  The two cocycles
differ by a coboundary $\delta \lambda$, and it is enough to prove
the result in the elementary case when the $0$-cochain $\lambda$ is
supported by one vertex $v_0$ of $T$. Again we have ``localized'' the
problem. Note that the transit-invariance for $1\to 4$ moves
establishes what we need in the special case when $v_0$ is the new
vertex after the move. Then we reduce the general case to this special
one, by means of the transit-invariance.

\noindent We may assume that $T$ is almost-regular.
As before, we probably have to generically
change the full cocycles in order to guarantee full transits. Again
this is not troublesome by the same continuity argument. 
Hence, it is enough to
show that, up to decorated moves, we can modify Star$(v_0,T)$ of a given
vertex of $T$ to reach the star-configuration of the special situation.
But Star$(v_0,T)$ is determined by the triangulation of its boundary, that
is of Link$(v_0,T)$, which is homeomorphic to $S^2$. So it is enough to 
control the link's modifications.
One sees that, by performing $2\to 3$ and $1\to 4$ moves
around $v_0$ in such a way that $v_0$ persits, their trace on    
Link$(v_0,T)$ are $1\to 1$ moves (2-dimensional analogues of the
$2\to 3$ moves) or $1\to 3$ moves (2-dimensional analogues of the
$1\to 4$ moves). It is well-known that these moves are sufficient to
connect any two triangulations of a given surface, so we are almost
done. We only have to take into account the technical complication
due to the fact that, in our situation,  
the  $1\to 4$ moves must be completed to moves of distiguished 
triangulations of $(W,L)$, that is we need to involve some edges
of $H$.

\medskip

\noindent The proof of Proposition \ref{propo1}, whence of the 
main Theorem \ref{teo1}, is now complete.

\medskip  
 
\subsection{Some properties of $K_N(W,L,\rho)$} \label{properties}

An essential question is the actual dependence of $K_N(W,L,\rho)$ on
$\rho$. Here are some very partial answers. The reader is advised to consult \S \ref{App} before reading this section.

\medskip

\noindent 
{\bf Projective invariance.} Take any $\mathcal{T}$ for $(W,L,\rho)$,
having as full co-cycle $z=\{z(e)\}$; $\{t(e)\}$ and $\{x(e)\}$ are as
before. For any $\lambda \neq 0$ we can turn $z$ into
$\{z_{\lambda}(e)\}$, where $\{t_{\lambda}(e)\}=\{t(e)\}$ and
$\{x_{\lambda}(e)\}=\{\lambda x(e)\}$.  In this way we get a
decoration $\mathcal{ T}_{\lambda}$ for some $(W,L,\rho_{\lambda})$.

\begin{prop} \label{invproj} For each $\lambda \neq 0$,
$$K_N(W,L,\rho_{\lambda})=K_N(W,L,\rho)\ .$$
\end{prop}

\noindent {\it Proof.} 
Since the functions $h$ and $\omega$ in the quantum dilogarithm
c-$6j$-symbols are homogeneous of degree $0$, each tetrahedron
contributes (via its c-$6j$-symbol) by $\lambda^{2p}$ to the state sum
$K_N(\mathcal{T}_{\lambda})$. Moreover, denoting by $r_i$ the number
of $i$-simplices of $T$ (remark that $r_0$ is also the number of edges
of $H$), the Euler characteristic of $W$
equals:
$$\chi(W) = r_0 - r_1 + r_2 - r_3 = 0\ .$$ 
As $r_2 = 2r_3$ we get $r_3 = r_1 - r_0$. In the formula of $K_N(\mathcal{T}_{\lambda})$, there are
only $r_1 - r_0$ edge contributions coming from $T \setminus H$, each
one being equal to $\lambda^{-2p}$. This concludes.

\medskip

\noindent{\bf Duality.} 
Let $\mathcal{ T} = (T,H,\mathcal{ D})$ and $z$ be as above. Let us
denote by $z^*$ the complex conjugate full co-cycle, $\rho^* =[z^*]$,
$\mathcal{ D}^* = (b,z^*,c)$, and $\widehat{W}$ the manifold with the
opposite orientation.

\begin{prop} $(K_N(W,L,\rho))^* = K_N(\widehat{W},L,\rho^*)$. 
Hence, if $\rho$ is ``real'', that is if it can be represented by
$B(2,\mr)$-valued cocycles, then
$$ (K_N(W,L,\rho))^*=K_N(\widehat{W} ,L,\rho)\ . $$
\noindent In particular
$$(K_N(W,L))^*=K_N(\widehat{W} ,L)\ . $$
\end{prop}

\noindent {\it Proof.} 
A change of orientation of $W$ turns the quantum dilogarithm
c-$6j$-symbol $t(\mathcal{ D},d)$ into $^T\bar{t}(\mathcal{ D},d)$. But
Proposition \ref{unitarite} shows that
$$^T\bar{t}(\mathcal{ D},\alpha,d) = 
\left( t(\mathcal{ D}^*,-\alpha,d) \right)^*.$$
\noindent Since the state sums $K_N$ are $N$'th-powers of 
weighted traces (see Lemma \ref{indbranch}), they do not depend on
$\alpha$. So, we get the conclusion.

\section{ The invariants $K_N(Z,S,\rho)$ and 
$K_N(Y,\phi,\rho)$ }\label{ZS}

Let $(Z,S,\rho)$ be as usual, and $W$ be 
the closed manifold obtained by
{\it Dehn filling} of each boundary component of $Z$ along $S$. If $L$
denotes the link of surgery cores, then $S$ becomes a family of
meridians of the $L_i$'s and $\rho$ is a flat bundle on $W\setminus L$
with, in general, a {\it non-trivial holonomy} along these
meridians. The aim of this section is to define state sum invariants
for the triple $(Z,S,\rho)$. Our construction consists in cutting up
pieces of any distinguished triangulation of $(W,L)$,
turning it into a decorated cell decomposition of $Z$ with
tetrahedron-like building blocks, which may be used to again define our
previous state sums.
\medskip

\begin{figure}[ht]

\begin{center}

\scalebox{1}{\input{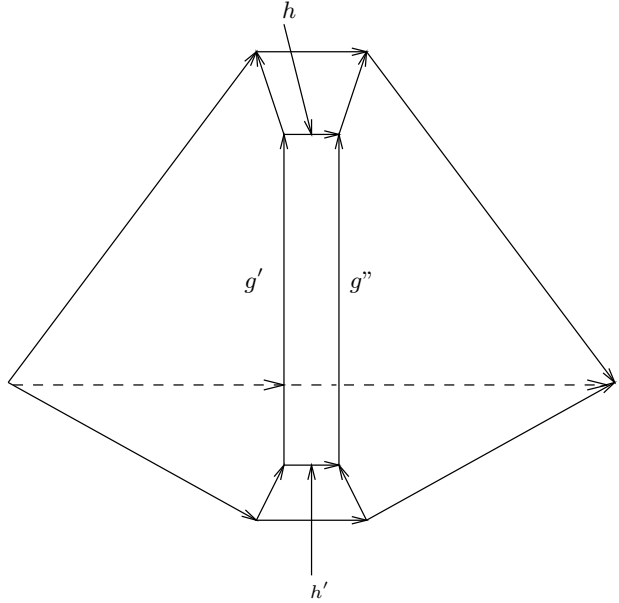}}

\end{center}

\caption{\label{decoup1} an elementary piece of $C$ }

\end{figure}

\noindent Let $(T,H)$ be a distinguished triangulation of $(W,L)$
 as above, with branching $b$ and charge $c$. For each edge $e \in H$
 consider an open neighborhood in each tetrahedron $\Delta \in T,e \in
 \Delta$, with the shape of a cylinder over a triangular basis. Let
 $U(H)$ be the union of all such polyhedra, for all edges. If $W' = W
 \setminus rD^3$, then $(T,H)$ induces a natural cell-decomposition
 $C$ of $Z = W' \setminus$ Int $U(H)$, with building blocks $\Delta'$
 as in Fig. \ref{decoup1}. Our first aim is to get a decoration of
 $C$.

\medskip

\noindent The edges of $T$, which subsist in $C$, are still $b$-oriented. In
each ``little'' face $l$ of $C$ coming from the link of the vertices
of $T$, the edges lying in the old faces of $T$ inherit their
$b$-orientation. Consider now the (new) rectangular face $f$ of some
$\Delta'$ at the place where an edge $e \in H$ has been
removed. Orient the two parallel copies of $e$ as $e$ itself. The two
other little arcs in $\partial f$ are in the boundary of some $l$:
give them the orientation opposite to that induced by the other ``long'' edges
of $l$. We thus have a canonical branching $b_C$ on $C$. Moreover, one
can represent $\rho$ by a \emph{full} $B$-$1$-cocycle $z_C$ on $C$, which
follows from the same arguments used in Corollary \ref{fullOK}. See
again Fig. \ref{decoup1}, where latin letters denote the values of
$z_C$ on $\partial f$.

\smallskip

\begin{figure}[ht]

\begin{center}

\scalebox{1}{\input{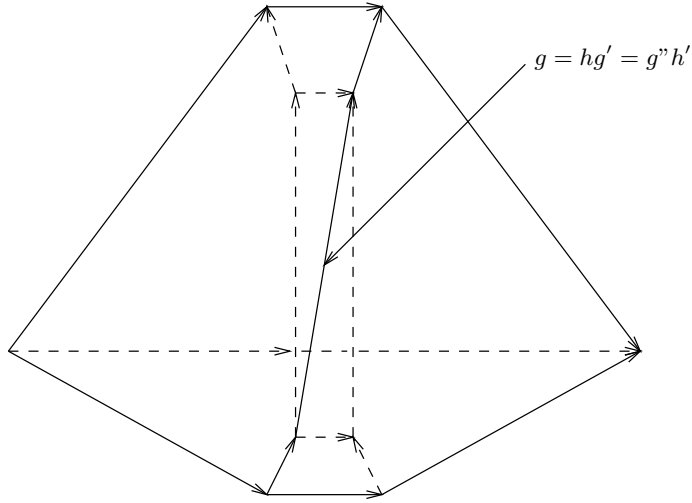}}

\end{center}

\caption{\label{decoup2} a model of $\widetilde{\Delta}$}

\end{figure}

\noindent For each 
$\Delta'$ and each rectangular face $f \in \Delta'$, there is
a unique diagonal arc in $f$ that joins the source and the sink of the
branching of $f$, and it is naturally oriented. Clearly, each
$\Delta'$ may then be seen as a deformation of a branched tetrahedron
$\widetilde{\Delta}$. Namely, the edges of $\widetilde{\Delta}$ are
obtained by forgetting some little edges and the parallel copies of
the edges of $H$ in $\Delta'$, all pictured in dashed lines in
Fig. \ref{decoup2}, and straightening the others. There is also a full
$B$-$1$-cocycle $z^{\tilde{\Delta}}$ on $\widetilde{\Delta}$ induced by
$z_C$. Remark that, in Fig. \ref{decoup2}, when the bundle extends to the
whole of $W$, then $h=h'=1$ and $g'=g''=g$. Each edge of
$\widetilde{\Delta}$ inherits the charge of the edge of $\Delta$ it
comes from. Finally, we thus have constructed a set
$\widetilde{\mathcal{T}}$ of decorated branched tetrahedra
$(\widetilde{\Delta},b^{\widetilde{\Delta}},z^{\widetilde{\Delta}},
c^{\widetilde{\Delta}})$
from $(T,H,b,c)$ and $z_C$, which is sufficient to define state sums
$K_N(\widetilde{\mathcal{T}},d)$ via formula (\ref{invariant})
below. Denote $\mathcal{D}^{\widetilde{\Delta}} =
(b^{\widetilde{\Delta}},z^{\widetilde{\Delta}},c^{\widetilde{\Delta}})$.

\medskip

\medskip

\noindent Note that one can glue the tetrahedra $\widetilde{\Delta}$ with the pairings of faces coming from $T$. It results in a compact polyhedron
$\widetilde{T}$ which is simplicially equivalent to $T$, with $L$
realized as an Hamiltonian subcomplex $\widetilde{H}$ inside. But $\tilde{z} = \{ z^{\Delta} \}$, which is defined separately on each $\widetilde{\Delta}$, do not have, in general, coinciding values along $\widetilde{H}$. We have to imagine that the edges of $\widetilde{H}$
initially defined a $1$-cycle in $\partial C = \partial Z$ (we probably
have to add some little edges from the faces $l$ to connect the
diagonals, but this does not matter), which is the union of broken
lines wound around $L$.

\noindent 
We say that $\widetilde{\mathcal{T}}$ is a decoration of $(Z,S,\rho)$.

\begin{teo}\label{t1prime} The state sum formula:
\begin{eqnarray} \label{invariant}
K_N(\widetilde{\mathcal{T}},d)=\left( N^{-r_0}\sum_{\alpha}
\prod_{\widetilde{\Delta}}
t^{\widetilde{\Delta}}(\mathcal{D}^{\widetilde{\Delta}},
\alpha,d)\prod_{e\in
E(\widetilde{T})\setminus E(\widetilde{H})} x(e)^{-2p/N} \right)^N\ .
\end{eqnarray}
\noindent defines a family of complex valued invariants 
$\{K_N(Z,S,\rho)\}$, $N$ being any odd positive integer. 
They recover the family $\{ K_N(W,L,\rho) \}$ when $Z=M$, 
each $s_i$ is a meridian of $L_i$, and $\rho$ extends through $S$.
\end{teo}

\noindent {\it Proof.} The only modification in the proof of Theorem \ref{teo1}
is that we have to verify that there is a well-defined state sum
transit-invariance for $K_N(\widetilde{\mathcal{ T}},d)$. But this is
a straightforward consequence of the fact that $2 \to 3$, $0 \to 2$
and $1 \to 4$ moves, which are purely local, apply without change on
$\widetilde{T}$.

\noindent 
Furthermore, when $Z=M$, each $s_i$ is a meridian of $L_i$ and $\rho$
extends through $S$, then $\tilde{z}$ is well-defined along
$\widetilde{H}$. This is the situation described in
fig. \ref{decoup2}, when $h=h'=1$.

\bigskip

\noindent Next we give the definition of the invariant $K_N(Y,\phi,\rho)$. Let
$Y$ be any oriented compact 3-manifold with non empty, non spherical and not necessarily
toral boundary. Set $\partial Y = \Sigma_1 \coprod \ldots \coprod
\Sigma_k$. The genus of $\Sigma_j$ is denoted by $g_j \geq 1$. For
each $g_j$, fix a base surface $\mathcal{H}_j$ regarded as the boundary of
a fixed handlebody of genus $g_j$, with a given complete system of 
meridians $\{ m_{ji}\}_{i=1,\dots ,g_j}$.

\begin{defi} 
\label{parametrisation} A {\rm symplectic parametrization} 
$\phi_j$ of $\Sigma_j$ is 
the equivalence class of an homeomorphism
$$p_j: \mathcal{H}_j \to \Sigma_j,$$ 
\noindent considered up to
self-homeomorphisms $\theta$ of $\mathcal{H}_j$ such that $\theta(m_{ji})
= m_{ji}$ (up to isotopy), and $\theta_*: H_1(\mathcal{H}_j;\mz) \to
H_1(\mathcal{H}_j;\mz)$ is an isometry for the standard symplectic form.
\end{defi}

\noindent A parametrization $\phi=(\phi_1,\ldots,\phi_k)$ of $\Sigma$ allows one to
define a triple $(Z,S,\rho_Z)$. Indeed, let $W$ be obtained by Dehn
filling of $Y$ along the system of curves $\phi_j(m_{ji})$, and $L$ be the surgery cores specified as the symplectic
duals to the $\phi_j(m_{ji})$'s. Then $Z = W \setminus$ Int $U(L)$, $S =
\{\phi_j(m_{ji})\}$ is a system of meridians of $L$ and $\rho_Z$ is the
restriction of $\rho$ to $Z$ (which embeds in $Y$). We thus have the
following generalization of Theorem \ref{t1prime}:

\begin{teo}\label{t2prime} 
One can define a family of complex valued invariants
$\{K_N(Y,\phi,\rho)\}$, $N$ being any odd positive integer, which may
be computed as a state sum from the triple $(Z,S,\rho_Z)$ constructed
above. Clearly, it specializes to the family $K_N(Z,S,\rho)$ when
$\partial Y$ is a union of tori.
\end{teo}

\section {Specializations of the invariants} \label{specialize} 

Purely formal considerations allow us to specialize
to closed manifolds or to {\it framed} links.
\smallskip

\noindent 
{\bf Closed manifolds invariants.}  Given $(W,\rho)$ as before, take
$L_0$ the {\it unknot} in a 3-ball in $W$. Clearly
$$K_N(W,\rho):=K_N(W,L_0,\rho)$$ 
\noindent 
is a well-defined invariant of $(W,\rho)$. Taking $\rho$ the trivial
flat bundle, we get a purely topological invariant of $W$.
\medskip

\noindent{\bf Framed links invariants.}  If $L$ is a link in $Y$, $X =
Y\setminus$ Int $U(L)$ and $S$ is made by longitudes of the components
of the link $L$ in $Y$, then we can interpret
$$K_N(L,S;Y,\phi,\rho):=K_N(X,(S,\phi),\rho)$$ 
\noindent as an invariant of the {\it framed} link $(L,S)$ in 
$(Y,\phi,\rho)$.
\medskip

\noindent More interesting specializations of the invariants arise by using 
special kinds of flat bundles $\rho$. We have already noticed in \S \ref{full}
that $B$ contains the abelian subgroups $Par(B)$, which is 
isomorphic to $(\mc,+)$, and $C(B)$, isomorphic to the multiplicative 
$\mc^*$. Moreover $U(1)$ naturally embeds into $C(B)$.

\smallskip

\noindent {\bf The Euler class.} Since $B$ retracts onto $U(1)$,
any principal $B$-bundle (forgetting the flat structure, if any) 
admits a group reduction to $U(1)$. Thus any $B$-bundle on $Y$ is
{\it topologically} classified by its {\it Euler class} in $H^2(Y;\mz)$.
Indeed, as $\pi_2(B)=0$, the obstruction to define
a global (continuous) section of it lies in $H^2(Y;\pi_1(B))\cong H^2(Y;\mz)$. 

\medskip

\noindent {\bf Induced $C(B)$-bundles and abelian-type invariants.} 
Since we are only interested in conjugation classes of representations in $B$, below we will often abuse of notations by writing fundamental groups of manifolds without any reference to basepoints. 

\noindent Any representation
$\rho : \pi_1(Y)\to B$ induces a representation $\rho':\pi_1(Y)\to C(B)$,
just by forgetting the upper-diagonal entries. Since $C(B)$ is Abelian,
$\rho'$ is trivial on the kernel of the natural projection
$p: \pi_1(Y)\to H_1(Y;\mz)$, that is $\rho'$ factorizes through $H_1$.
Thus, in a sense, the group $H^1(Y;C(B))\cong H^1(Y;\mc^*)$ 
is the ``core'' of our discussion, and there is a natural map
$$\mu: H^1(Y;C(B))\to H^1(Y;B)\ .$$ 

\noindent We have already noticed in \S \ref{full}
that $H^1(Y;Par(B))\cong H^1(Y;\mc)$, and that it maps via 
$$\alpha: H^1(Y;Par(B))\to H^1(Y;B)\ .$$ 

\noindent We say that any invariant of the form
$K_N(Y,\phi,\rho=\mu(x))$ (resp. $K_N(Y,\phi,\rho=\alpha(x))$) is an 
invariant of {\it multiplicative} (resp. {\it additive}) abelian type.
\medskip

\noindent {\bf Abelian-type invariants and the hyperbolic 
 $VCS$-conjecture.} Assume that $(W,L)$ has some
``hyperbolization'' $\mathcal{F}$ in the sense we have already
stipulated. In order to make more consistent the {\it Qualitative
$VCS$-conjecture} stated in the introduction, it would be very useful
to determine flat $B$-bundles intrinsically related to the hyperbolic
manifold $\cal F$. To make this claim precise, we are going to elaborate it
in the special case of a hyperbolic knot $L$ in $S^3$.  So, either
$\mathcal{F} = S^3\setminus L$, or $\mathcal{F}$ is the (compact)
result of a hyperbolic surgery along $L$ and $\mathcal{G}$ denotes the
geodesic core of the surgery. Using the notations of \S
\ref{charges}, the (refined) scissors congruence class
$\widehat{\beta}(\mathcal{F})$ and $\widehat{\beta}(\mathcal{F},\mathcal{G})$ 
are intimately related to the
{\it affine structure} induced by $\mathcal{F}$ on the boundary torus
$\partial U(L)$. This structure is Euclidean exactly in correspondence
with the {\it complete} hyperbolic structure on $S^3\setminus L$, and the
linear part of its holonomy is then the identity. We interpret this fact as
the main source of the role of the trivial flat bundle $\rho$ in the
current volume conjecture for hyperbolic knots in $S^3$. In the case of a hyperbolic 
surgery, we are considering a suitable {\it non-complete} hyperbolic structure
on $S^3\setminus L$ (close to the complete one) and $\mathcal{F}$ is
just given by its completion.  The linear part of the holonomy of the
affine structure is now a non-trivial representation $h:
\pi_1(\partial U(L))\to \mc^* \cong C(B)$.  In particular, $h$ is non
trivial on the canonical longitude of $L$ (which is homologically
trivial), so it cannot be extended to the whole of $\pi_1(\mathcal
{F})$. Nevertheless, if $m$ denotes (the copy in $\cal F$ of) a
meridian of $L$, then $h$ extends to $h: \pi_1(\mathcal{F} \setminus m)
\to C(B)$; in fact $S^3 \setminus$ Int $U(L\cup m)$ is a
$\mz$-homology-$\partial U(L) \times [0,1]$, and $h$ is trivial on the
essential curve on $\partial U(L)$ which is filled by a 2-disk in
$\cal F$, via the surgery. This construction is completely
canonical. One could interpret $m$ as a sort of {\it cut locus} for
the original representation $h$. We dispose now of specific
abelian-type invariants for $\mathcal{F} \setminus$ Int $U(m)$ (or
$\mathcal{F} \setminus$ Int $U(m \cup \mathcal{G})$) which are canonically associated to
$(\mathcal{F},\mathcal{G})$, and we can use them to make the
$VCS$-conjecture precise in this case.
A suitable generalization of this
``cut-locus theory'' for the affine holonomy on $\partial U(L)$ could
also play an important role in specifying a general form of the $VCS$-conjecture. We shall
develop these considerations in \S \ref{GR}.
\medskip

\noindent {\bf (Singular Homology)-derived invariants.} Concrete examples of
abelian-type invariants can be derived from the ordinary singular 
cohomology of $Y$. 

\noindent Fix $\lambda \in \mc$. By means of the
exponential map $\exp_{\lambda} : \mc \to \mc^*$, 
$\exp_{\lambda}(z) = \exp(\lambda z)$, 
we get a map $\exp_{\lambda}: H^1(Y;\mc) \to H^1(Y;C(B))$ (which restricts
to the free abelian group $H^1(Y;\mz)/Tors$).
 
\noindent So the ordinary cohomology group $H^1(Y;\mc)$ maps in two ways 
to 
$H^1(Y;B)$, called respectively {\it multiplicative} 
(just take  $\mu_{\lambda} = \mu \circ \exp_{\lambda}$), and 
{\it additive} (take the above $\alpha$).
 
\noindent For $\lambda = h = 2i\pi/k$, $k\in \mn^*$, one also has the 
embedding
$\exp_h : \mz/k\mz \to U(1)$, $\exp_h([m]) = \exp_h(m)$, hence a map 
$\exp_h: H^1(Y;\mz/k\mz) \to H^1(Y;U(1))$, and finally a map also denoted
by $\mu_h: H^1(Y;\mz/k\mz) \to H^1(Y;B)$. Clearly the two $\mu_h$'s ``coincide''
on $H^1(Y;\mz)/Tors$ via the natural map 
$H^1(Y;\mz)/Tors \to H^1(Y;\mz/k\mz)$. 

\begin{remark} {\rm For $x\in H^1(Y;\mc)$, $\mu_{\lambda}(x)$ and 
$\alpha(x)$ should lead to different invariants. For example, let $L$ be a 
knot in $S^3$, $S$ its canonical longitude, $x$ the Poincar\'e dual
of a generator of $H_2(Z,\partial Z;\mz)\cong \mz$, $Z= S^3 \setminus$ Int $U(L)$.
Then, for any $a\in \mc^*$, 
$K_N(Y,\phi,\alpha(ax))= K_N(Y,\phi,\alpha(x))$; this is a consequence
of the ``projective invariance'' \ref{invproj}. On the other hand, 
$K_N(Y,\phi,\mu_{\lambda}(ax))$ depends on $a$. Apparently the 
multiplicative approach is more sensitive.}
\end{remark}

\begin{proble} {\rm Given $x\in H^1(Y;\mc)$, carefully analyze the 
dependence of $K_N(Y,\phi,\mu_{\lambda}(ax))$ on $(\lambda,a)$.}
\end{proble} 
  
\noindent{\bf Seifert-type invariants.} Here we specialize the last 
construction to the case of an 
{\it oriented} link $L = \coprod_{i=1}^n L_i$ in a $\mz$-homology
sphere $\mathcal{H}$. Let $F_i$ be any oriented Seifert surface for
the component $L_i$ in $Z_i=\mathcal{H} \setminus$ Int $U(L_i)$, $Z =
\mathcal{H} \setminus$ Int $U(L)$, and $l = (l_1,\ldots,l_n)$, with
$l_i$ a longitude for $L_i$.  As usual, $l$ represents a {\it framing}
of $L$ and it can be encoded as an element of $\mz^n$. The surface $F_i$
represents a class in $H_2(Z_i,\partial Z_i;\mz)$ which is the
Poincar\'e dual of a class $f'_i\in H^1(Z_i;\mz)$; by the embedding
$j_i$ of $Z$ in $Z_i$, we get $f_i := j_i^*(f'_i) \in
H^1(Z;\mz)$. Finally, for any $c=(c_1,\ldots,c_n)\in \mz^n$, define
$$x(c) = \sum_{i=1}^n c_i [f_i] \in H^1(Z;\mz).$$
\noindent Set $\lambda = 1$ or  $\lambda = h = 2i\pi/k$, $k\in \mn^*$. 
Any invariant of the form $K_N(Z,l,\mu_{\lambda}(x(c)))$ or 
$K_N(Z,l,\alpha(x(c)))$ is called a {\it Seifert-type invariant}.

\smallskip

\noindent We believe that Seifert-type invariants (when $\mathcal{H} = S^3$)
are good candidates in order to reconduct the full set of coloured
Jones polynomials in the main stream of the quantum hyperbolic
invariants and produce new interpretations for them. Our work on this
matter is in progress (see \cite{12}). Here we limit ourselves to some
considerations in that direction.

\smallskip

\noindent 
Set $\mathcal{H} = S^3$, $L$ a link in $S^3$, and 
$Z= S^3\setminus$
Int $U(L)$ as usual. Let $J(L,(c_1,\ldots,c_n))_{\vert t}$ be the
non-framed version of the coloured Jones polynomial of the link $L$ in $S^3$
for $sl(2,\mc)$, evaluated in $t \in \mc^*$ \cite{29}. Then, each
component $L_i$ is framed by its \emph{canonical} longitude and
coloured with the $c_i$-dimensional irreducible representation of type
$1$ of the simply-connected restricted integral form
$U_q^{res}(sl(2,\mc))$ \cite[\S 9.3]{30}, specialized in $q=t^{\frac{1}{2}}$ a fixed square root of $t$; 
$J(L,(c_1,\ldots,c_n)) \in \mz [t^{\pm \frac{1}{4}}]$. 
Moreover, $J(L,(c_1,\ldots,c_n))$ is normalized by the value
on the unknot 
$J(O,c)_{\vert t} = [c]_{t}= \frac{t^{c/2} - t^{-c/2}}{t^{1/2} -
t^{-1/2}}, t \ne 1$. Define the rational function:
$$J^{\circ}(L,(c_1,\ldots,c_n)) = \frac{J(L,(c_1,\ldots,c_n))}{[c_1]},$$
so that, in particular, $J^{\circ}(O,N)_{\vert \omega} = 1$.

\noindent  Recall that in \cite{2}, 
Kashaev presented an invariant of links in $S^3$, well-defined up to $N$'th roots of unity, using an $R$-matrix
derived from the canonical element of the Heisenberg double of the
Weyl algebra, which is a (twisted) quantum dilogarithm
\cite{4,16}. He also claimed that the $N$'th power of this invariant should coincide with
the specialization for links in $S^3$ of his previous state-sum
proposal \cite{1}.  Murakami-Murakami \cite{8} showed that this
$R$-matrix can be made into an enhanced Yang-Baxter operator
equivalent to the one of the coloured Jones polynomial
$J^{\circ}(L,(N,\ldots,N))$ evaluated in $\omega$. Then, combining both facts in our context would give an equality $K_N(Z,l_0,\mu_{\lambda} (x(0))) = (J^{\circ}(L,(N,\ldots,N))_{\vert \omega})^N$, where $l_0$ denotes the system of canonical longitudes of 
$L$ in $S^3$. But there is at present no proof of the claimed coincidence between
the ``$R$-matrix'' invariants and the ``state-sum'' ones.  In fact, we are able to obtain this coloured Jones polynomial with a slight modification of our setup, but in order to get an actual coincidence with the Seifert-type invariants, we believe that some computations need to be better understood. This fact deserves a careful analysis, postponed to
\cite{12}, with the hope that it could reveal itself as the prototype
for a full generalization.  Note that $K_N(Z,l,\alpha(x(c)))$ should be
sensitive to the orientation of $L$ when $i \geq 2$, since we can use
Proposition \ref{invproj} only for simultaneous changes of the
orientations of all the link's components. The evident \emph{cabling}
information contained in the ``colours'' $c_i$ of $x(c)$ makes the
Seifert-type invariants sensitive to the linking matrix of $L$. Also,
Seifert-type invariants $K_N(Z,l,\mu_h(x(c)))$, $h\neq 1$, yield a
``Symmetry Principle'' for colours, in the style of \cite{29},
\cite[\S 4]{31}. The following rather informal problem contains, 
nevertheless, an implicit conjecture about the expected result:

\begin{proble} \label{q1} 
{\rm Express the framed coloured Jones polynomial 
$J(L,l,(c_1,\ldots,c_n))$ in terms of the 
Seifert-type invariants $K_N(Z,l,\mu_h x(c_1,\ldots,c_n))$.}  
\end{proble}

\begin{remark} \label{U(1)} {\rm  Since $B$ retracts on $U(1)$, a 
solution of Problem \ref{q1} would also establish, in the case of 
\emph{topologically trivial} B-bundles, a correspondence between 
our results and Rozansky's work \cite{32} on $U(1)$-reducible 
SU(2)-connections over complements of links in $S^3$ and the 
Jones polynomials.}
\end{remark}

\section{Speculations on the $VCS$-conjecture and on (2+1) quantum gravity }
\label{GR}

\subsection{More VCS-Conjecture}\label{VCS}
\noindent Let us adopt the notations of \S \ref{charges}.
Let $\cal F$ be an oriented finite volume  hyperbolic 3-manifold. 
A fundamental result of \cite{13,14} and \cite{15} relates the (refined) scissors congruence class of $\cal F$ and its
$VCS$-invariant.

\begin{teo}\label{*} $$R(\widehat{\beta} ({\cal F}))= i 
\left( Vol({\cal F})+iCS({\cal F}) \right)\ .$$
\end{teo}

\noindent If $\cal F$ is non-compact we know that $\widehat{\beta} ({\cal F})$
can be explicitely expressed as $\widehat{\beta} (T)$, using any
geometric ideal triangulation $T$ of $\cal F$. When $\cal F$ is
compact it is subtler to express the scissors class by means
of ``ideal triangulations'' of $\cal F$: even the meaning of such a notion
is not evident a priori. However, if $\cal F$ is the result of 
a hyperbolic Dehn surgery with geodesic core $\cal G$, one has, in a natural
way, ideal triangulations $T$ of $\mathcal{F}\setminus \mathcal{G}$ which lead
to $\widehat{\beta} (\mathcal {F},\mathcal{G})$ (see \S \ref{charges}).
Using it and \cite{23,15}, we have:

\begin{teo}\label{**}
$$ R(\widehat{\beta} (\mathcal {F},\mathcal{G}))
= i VCS(\mathcal{F}) + \frac{i\pi}{2}\sum_{i=1}^n
\lambda(L_i) \in \mc/\pi^2\mz\ := iVCS(\mathcal{F},\mathcal{G})\ .$$
\noindent  where $\lambda(L_i)$ denotes the
length of the closed geodesic $L_i$ plus $i$ times its torsion.
\end{teo}

\smallskip

\noindent Consider now a hyperbolic knot $L$ in $S^3$. Let $\cal F$ be
either $S^3\setminus L$ or the result of a hyperbolic Dehn surgery along
$L$ with geodesic core $\cal G$. In \S \ref{specialize} we have shown
how to associate to $\cal F$ abelian-type invariants which we denote here 
$K_N(\cal F)$ or 
$K_N(\mathcal {F}, \mathcal{G})$ respectively. If $\cal F$ is non
compact $K_N(\mathcal {F})$ is just the original Kashaev's $K_N(S^3,L)$. 
In the other case 
$K_N(\mathcal{F}, \mathcal{G}):=K_N(\mathcal{F} \setminus \ {\rm Int} \ U(m),m',\rho)$,
where $m$ is (a copy in $\cal F$ of) a meridian of $L$, $m'$ is a meridian of $m$, and the flat
bundle $\rho$ on $\mathcal{F} \setminus$ Int $U(m)$ is canonically associated to the affine structure induced
by the hyperbolic structure of $\cal F$ on $\partial U(L)$.
\noindent Moreover, the $N$-dimensional quantum dilogarithm 6j-symbols and the 
pentagon identities that they satisfy may be respectively interpreted as
deformations of the Euler dilogarithm and the five-term relation
(\ref{five-term}) of the classical Rogers dilogarithm. Indeed, both
may be recovered from them in the limit $N \to \infty$
\cite{4,5,16}; in the pentagon identities, this is the contribution
of the powers of $\omega$ which turn the Euler dilogarithm into the
Rogers one.
\noindent Hence we can make the $VCS$-conjecture precise at least for hyperbolic knots in
the 3-sphere, which looks in this case to be {\it morally
tautological}.

\begin{conge} \label{cknots} (VCS-conjecture for hyperbolic knots in $S^3$.) 

\noindent (a) When $\mathcal{F} = S^3\setminus L$ then
$$VCS(\mathcal{ F})= lim_{N\to \infty} (2\pi/N)\log (K_N(\mathcal{F}))\ .$$

\noindent (b) When $\mathcal{F}$ is compact
$$VCS(\mathcal{ F},\mathcal{G})= 
lim_{N\to \infty} (2\pi/N)\log (K_N(\mathcal{F},\mathcal{G}))\ .$$
\end{conge}

  \begin{remark} {\rm

(1) Recall that there are also non-complete hyperbolic structures on
$S^3\setminus L$ which lead, by completion, to $\cal F$ with a
{\it conical singularity} along $\mathcal{G}$. The above discussion applies
also to this situation.  
\smallskip

(2) When $\mathcal{F}=S^3\setminus L$ the conjecture fits
well with Theorems \ref{*} (non compact case). One would
like to extend it to any $\mathcal {F}= W\setminus L$ where $L$ is  
a  hyperbolic link in any 3-manifold $W$. But we prefer to 
prudently stay with a knot in $S^3$ because, in the general case, 
there is not a unique way to
extend the trivial flat $B$-bundle on $\partial U(L)$ to $W$; in the
spirit of the ``path integral'' ideology which we mention below, this could lead to non-trivial asymptotic contributions
different from the $K_N(W,L)$ one. 
\smallskip

(3) When $\mathcal{F}$ is compact
Theorem \ref{**}  seems to give the most appropriate ``classical'' counterpart,
because the quantum hyperbolic invariants should depend
also on $\cal G$, not only on $\cal F$.
\smallskip

(4) We believe to be still far from an equally effective formulation
of the $VCS$-conjecture for arbitrary hyperbolic links $L$ in an
arbitrary $W$. In particular we should have to develop a proper general
notion of ``representation cut locus'' (see \S \ref{specialize}). 
\smallskip

(5) The Conjecture \ref{cknots} is supported by some numerical 
computations, obtained by
using a somewhat formal application of the stationary phase method on
an integral formula for the coloured Jones polynomials. It is not
difficult to generalize such a formula for our state sum invariants,
but we stress that, in order to get an {\it actual} verification of the
conjecture in this way,  there are still many analytical problems related
to the study of its asymptoptic behaviour, even for elementary
examples of knots in $S^3$.}
\end{remark}

\subsection{(2+1) quantum-gravity}

\noindent 
Witten had already suggested a form of the VCS-Conjecture, in the
context of the matter-free Euclidean (i.e. Riemannian) continuation of
quantum gravity in $(2+1)$ dimensions, with negative cosmological
constant $\Lambda$ \cite{10}. This and other facts suggest that this theory is
a very meaningful heuristic background for the interpretation of our
invariants.
\medskip

\noindent 
Let us briefly recall how this goes. The equation of motions of
(matter-free) classical gravity with cosmological constant $\Lambda$
are obtained from the Euler-Lagrange equations associated to the
Einstein-Hilbert action $I_{EH}(\Lambda)$, varying the metric of the
underlying space-time manifold $W$. They say that the latter is
locally homogeneous, with a scalar curvature $K$ proportional to
$\Lambda$. With a view towards geometric applications, let us suppose
that $W$ is oriented, closed (for simplicity) and endowed with a Riemannian structure. Then, 
when $\Lambda < 0$, up to a constant conformal factor, $W$
is a \emph{classical} solution of $3$-dimensional gravity if it is a
complete hyperbolic manifold: $W \cong \mathbb{H}^3/\Gamma,$ where
$\Gamma$ is a discrete subgoup of $PSL(2,\mc)\cong SO(3,1)$.

\noindent 
Now, consider a principal bundle $P: SL(2,\mc) \rightarrow E
\rightarrow W$, with a connection $A \in sl(2,\mc) \otimes \Omega^1(E)$
over it. Following Chern-Simons \cite{33}, for any global continuous section s of $E$ (which exists since $\pi_1(SL(2,\mc)) = \pi_2(SL(2,\mc)) = 0$ and $SL(2,\mc)$ is connected), define the
\emph{Chern-Simons $3$-form} for $(A,s)$ as:

$$ CS(A,s) = \frac{-1}{8\pi^2} s^* \left( tr(A \wedge dA + \frac{2}{3} A \wedge A \wedge A) \right) \in \Omega^3(W,\mc),$$
where $tr$ denotes any non-degenerate invariant bilinear form on $sl(2,\mc)$. Then the so-called \emph{Chern-Simons action}:

\begin{eqnarray} 
 I_{CS}(A) := I_{CS}(A,s) = \frac{-1}{8\pi^2} \int_W s^* \left( tr(A \wedge dA + \frac{2}{3} A
 \wedge A \wedge A) \right) \in \mc/\mz,\nonumber
\end{eqnarray}
is well-defined mod $\mz$, for $I_{CS}(A,s)$ only depends on $s$ via the \emph{degree} of the induced map $\bar{s}: W \rightarrow SL(2,\mc)$. Set $\Lambda = -1$. This gives no
restriction for the above classical solutions, since any discrete and
faithful representation of $\pi_1(W)$ lifts to $SL(2,\mc)$. In
\cite{10}, it is shown that $I_{EH} := I_{EH}(-1)$ is proportional to
an action of Chern-Simons type in two essentially \emph{different} ways,
which we denote by $I$ and $\tilde{I}$. Each may be constructed by
writing $A$ in terms of a frame field $e$ of $TW$ and an
$SL(2,\mr)$-connection $w$ over $E$, and using respectively the
Killing form and the other non-degenerate invariant form of the Lie
algebra $sl(2,\mc)$.

\smallskip

\noindent 
Then, Witten's prediction holds for the Feynman integral over the
field variables, with Lagrangian $\hat{I}_{\hbar} = \frac{1}{\hbar}(I
+i\tilde{I})$ (where $\hbar$ is real). 
More precisely, for a closed and complete hyperbolic manifold $W$, the
small-$\hbar$ limit of the partition function

$$Z(W) = \int \exp (- \hat{I}_{\hbar} ) \ \mathcal{D}e \ \mathcal{D}w$$

\noindent should have an exponential growth rate equal to $VCS(W)$. 
For an arbitrary (closed) $W$, one expects $Z(W)$ to be localized around the
stationary solutions of the variation of $\hat{I}$ with respect to
$e$ and $w$. Moreover, by \cite{10}, these solutions should be the moduli space
$\chi_{SL(2,\mc)}$ of flat $SL(2,\mc)$-connections over $W$. Then,
from results of Dupont \cite{34} relating the Chern-Simons classes of flat connections to
volumes of representations, this prediction is also definitely well shaped for arbitrary closed
$3$-manifolds.

\medskip

\noindent 
The Turaev-Viro state sum invariants, based on the quantum group
$U_q(sl(2))$ with $q=\exp(2\pi i/k)$, $k$ a positive integer, are
commonly regarded as $q$-deformations of the Euclidean continuation of
(2+1) quantum gravity with {\it positive} $\Lambda$ and gauge group
$SU(2)$ \cite[\S 11.2]{35}. They are partition functions of a well-defined 
TQFT; in this context,
$\Lambda$ is ``quantized'' and it is related to $q$ by $2\pi/k =
\sqrt{\Lambda}$. So, the limit for $k\to \infty$ corresponds to
$\Lambda \to 0$.

\noindent 
Arguing similarly with respect to the action $I+ i\tilde{I}$, one may
look at the state sums invariants $K_N$, established in the present
paper, as parts of $q$-deformations of the Riemannian (2+1) quantum
gravity with {\it negative} $\Lambda$ and gauge group
$SL(2,\mc)$. Here, the relevant quantum group is the Borel subalgebra
$\mathcal{W}_N = BU_q(sl(2,\mc))$, with $q=\exp(2\pi i/N)$, and it takes
part only via its cyclic representation theory. Then $2\pi /N =
\sqrt{-\Lambda}$, and the limit for $N\to \infty$ should correspond to
$\Lambda \to 0$. From Conjecture \ref{cknots} 
(and its hypothetical generalizations) we expect that the $K_N$ invariants
capture some dominant part of the \emph{integrand} of the expected
\emph{full} partition functions. That is we should have a dominant
``rough'' partition function, for which one integrates over some
collection  $\chi_0$ of
irreducible components of Hom$(\pi_1,B)/B$, instead of the whole moduli space
Hom$(\pi_1,SL(2,\mc))/SL(2,\mc)$. Then, taking
inspiration from \cite{36}, we may think for instance of partition functions for pairs
$(Z,S)$, which formally look like:

$$F_N(Z,S) = \int_{\chi_0} K_N(Z,S,\rho)\ \mathcal{D}\rho\ .$$

\noindent Considering the closed $3$-manifold $W$ obtained by Dehn filling of $Z$ along $S$
with cores $L$, the link $L$ may be considered as ``Wilson lines'' for the integrand. 
In order to get complete regularizations of these partition functions for some TQFT  
(similarly to the Turaev-Viro state sums), one could  try to restrict to some
(finite) particular subsets of Hom$(\pi_1,B)/B$, and consider
weighted sums of the $K_N$'s over them. A solution to Problem \ref{q1}  
would allow to recover in this way, for example, just the 
Reshetikhin-Turaev and the Turaev-Viro invariants for $sl(2,\mc)$.
\medskip

\noindent 
By the proof of the Melvin-Morton conjecture \cite{37,38} and
Problem \ref{q1}, it could be possible that the {\it absolute torsion}
$\tau(ad\circ \rho)$ (see \cite{17}) of the flat vector bundle on $W$
obtained by means of the adjoint representation of $B$ can be
recovered from the $K_N$ invariants in the limit for $N$ large. This is in agreement
with our proposal for $F_N$ and some speculations in \cite{36}. At
this point, let us recall that
$$\tau(ad\circ \rho)= \tau(ad\circ \rho,\xi_c),$$
\noindent where $\xi_c$ is the {\it canonical} Euler structure (i.e.
the Spin$^{\mc}$ structure) on $W$ with null Euler class
$e(\xi_c)=0\in H^2(W;\mz)$. Staying at a formal level, we would like
to stress that the same technology used to define $ K_N(W,L,\rho)$ can
be used to treat the absolute torsions. Many basic extra-structures on
3-manifolds such as {\it framings}, {\it spin structures} and {\it
Euler structures}, can be encoded in a very effective way by means of
branched spines \cite{20},\cite{25,39}.  In order to obtain also an
encoding of Euler structures in terms of the canonical {\it Euler
chains} carried by branched spines (as in \cite{39}), it is enough to
require that the distinguished decorated triangulations
$(T,H,\mathcal{ D})$ of $(W,L)$ satisfy the further condition:

\medskip

\noindent {\it For each vertex $v_0$ of $T$ consider the ``black'' region
$B(v_0)$ of $S^2(v_0)$, the 2-sphere around $v_0$, 
where the natural field on $W'_r$,
normal  to the branched spine $P(T)$, points inward the ball around $v_0$. 
Then
$$\chi(B(v_0))=0 \ .$$}

\noindent This condition ensures that the normal field
can be extended to a {\it non-singular field} on the whole of $W$, whence it
determines an Euler structure $\xi$ on $W$. The Euler class $e(\xi)$ of this
Euler structure can be explicitely represented by a canonical
2-co-chain on $P(T)$, and we know exactly how to modify $P(T)$ (hence
$\xi$)
by means of $2\to 3$ (branched) moves in order to make $e(\xi)=0$
(see \cite{20} for the details).
\medskip

\section{Quantum data} \label{App}

In this Appendix we state 
the algebraic results needed for the construction of our 
invariants. Most of them were already announced by R. Kashaev in 
\cite{1,3}. The proofs are often involved formal computations. Nevertheless,  they do not add anything to the understanding of the rest of this paper, so we refer to \cite{16} for the details. 
Recall that $\omega = \exp(2i\pi/N)$ for an odd positive 
integer $N$, and $N= 2p+1,\ p \in \mathbb{N}$. 
Fix the determination $\omega^{1/2} = \omega^{p+1}$ for its 
square root. We shall henceforth denote 
$1/2 := p + 1 \ \rm{mod} \ N$. All other notations for 
manifolds, triangulations and spines are as in the rest of 
the paper.

\medskip

\noindent Since the work of \cite{6}, we know that the quantum 
$6j$-symbols and, more generally, the associators of 
braided, spherical or balanced semisimple categories provide 
via the Pentagon equation a natural framework for 
constructing invariants of compact $3$-manifolds 
\cite{7,40,41}. In each of these cases, 
the construction heavily uses some symmetry properties of 
the associators for the group $\mathbf{S}_4$ of permutations 
over $4$ elements, which is the group of symmetries of a 
tetrahedron. But in our present situation, the number of 
irreducible cyclic representations of the Weyl algebra 
$\mathcal{W}_N$, in any fixed dimension, is infinite 
(see Proposition \ref{prop1}), and the $6j$-symbols 
$R(p,q,r)$ obtained in Proposition \ref{6j} 
(and then, any $6j$-symbol for the regular sequences of 
cyclic representations of $\mathcal{W}_N$) are not symmetric 
under the action of $\mathbf{S}_4$.\hfill\break

\noindent One could try to overcome this difficulty by considering 
``branched'' versions of the pentagon equation and possibly by
restricting to some topologically significant subset of 
branched moves such as the {\it sliding moves} (see \S \ref{branchings}); 
unfortunately, 
this does not work. Hence we are forced to find a 
generalization of $R(p,q,r)$ which is covariant under 
$\mathbf{S}_4$. There is a representation-theoretic 
interpretation of this extension, which is related to the 
action of $\mathbf{SL}(2,\mathbf{Z})$ on $\mathcal{W}_N$ (see \cite{16}, 
the discussion before Proposition \ref{6j} and Remark \ref{matrices}).

\smallskip

\noindent We first briefly recall some fundamental results 
concerning the \emph{Weyl algebra} $\mathcal{W}_N$. This is 
the unital algebra over $\mathbb{C}$ generated by elements $E,\ E^{-1},D$ 
satisfying the commutation relation $ED = \omega DE$.
It is well-known that $\mathcal{W}_N$ can be endowed with a 
structure of Hopf algebra isomorphic to the simply-connected 
(non-restricted) integral form of the Borel subalgebra of $U_q(sl(2,\mathbf{C}))$ \cite[\S 9-11]{30}, specialized in $\omega$, with the following 
co-multiplication, co-unity and antipode maps :
$$\begin{array}{l}\Delta(E)  = E \otimes E,\ \Delta(D) =  
E \otimes D + D \otimes 1,\\ \epsilon(E)=1,\ \epsilon(D)=0,\ S(E)  = E^{-1},\ S(D) = -DE^{-1}.\end{array} \nonumber$$
\noindent Consider now the injective representations of 
$\mathcal{W}_N$ on finite dimensional complex vector spaces, 
which are also called the \emph{cyclic representations} of 
$\mathcal{W}_N$ in the theory of quantum groups. 
Let us fix the dimension equal to $N$, and let 
$X$ and $Z$ be the $N \times N$ matrices with components: 
$$X_{ij} = \delta_{i,j+1},\  Z_{ij} = 
\omega^i \delta_{i,j},\label{matrices}$$
\noindent where $\delta_{i,j}$ denotes the symbol of 
Kronecker. Following \cite{3}, we define a 
\emph{standard $N$-dimensional representation} $p$ of 
$\mathcal{W}_N$ by :
$$ p(E) = t_p^{2}Z,\ p(D) = t_px_pX,\ t_p,x_p \in 
\mathbb{C}^*,$$
\noindent and we denote by $p^*$ (the \emph{complex conjugate} 
representation) and by 
$\bar{p}$ (the \emph{inverse} representation) the 
standard representations of $\mathcal{W}_N$ obtained by 
setting
\begin{eqnarray} t_{\bar{p}} = 1/t_p,\  x_{\bar{p}} 
= -x_p \nonumber \\ t_{p^*} = (t_p)^*, \ x_{p^*} = (x_p)^*.    
\nonumber
\end{eqnarray}
\noindent We immediately see from the definition of the comultiplication 
$\Delta$ that for any pair $(p,q)$ of standard representations, 
we have
\begin{eqnarray}
 p \otimes q(E^N) = (t_pt_q)^{2N} \ id \otimes id, \nonumber \\
 p \otimes q(D^N) = (t_pt_q)^N \left(t_p^Nx_q^N + 
\frac{x_p^N}{t_q^N}\right) \ id \otimes id,\label{form2}
\end{eqnarray}
\noindent where $id$ denotes the identity matrix in dimension $N$. 
In particular, for any standard representation $p$, the representation 
$p \otimes \bar{p}$ is no longer injective. Then we call \emph{regular} an 
ordered sequence of cyclic representations $(p_1,\ldots,p_n)$ of 
$\mathcal{W}_N$ for which any tensor product representation 
$p_i \otimes \ldots \otimes p_{i+j},\ 1 \leq i \leq n,\  1 \leq j \leq n 
-i$ is cyclic. The following proposition shows that the isomorphism classes of cyclic representations of 
$\mathcal{W}_N$ may be described using only standard representations. Moreover,
it implies that they are infinite in number and that their tensor product is not commutative. Hence they do not form either a semi-simple or a monoidal category.

\begin{prop} \label{prop1} i) The standard representations 
are irreducible and cyclic, and two standard $N$-dimensional 
representations $p$ and $q$ are equivalent if and only if

$$t_p^{2N} = t_q^{2N},\ t_p^Nx_p^N = t_q^Nx_q^N.$$

\noindent ii) Any irreducible cyclic representation of 
$\mathcal{W}_N$ is equivalent to a standard one. 

\noindent iii) Fix a branch of the $N$'th root. 
If $(p,q)$ is a regular pair of $N$-dimensional irreducible 
cyclic representations of $\mathcal{W}_N$, the representation 
$p \otimes q$ is equivalent to the direct sum of $N$ copies 
of the standard representation $pq$ defined by : 
\begin{eqnarray} t_{pq} = t_pt_q,\ x_{pq} = \left(t_p^Nx_q^N 
+ \frac{x_p^N}{t_q^N}\right)^{1/N}. \nonumber \end{eqnarray}
\noindent Such a representation is called a \emph{product} 
standard representation. 
\end{prop}
\noindent Define the space of \emph{intertwiners} for 
representations $\rho$ and $\mu$ of $\mathcal{W}_N$ acting 
respectively on the vector spaces $V_\rho$ and $V_\mu$ as
$$H_{\rho,\mu} = \{ U : V_\rho \rightarrow  V_\mu \vert \ 
U\rho(a) = \mu(a)U,\ a \in \mathcal{W}_N \}.$$
\noindent Then, the claims $ii)$ and $iii)$ of Proposition 
\ref{prop1} imply that for any irreducible cyclic 
representations $p,q,r$ of $\mathcal{W}_N$ with $(p,q)$ a 
regular pair, the dimension of $H_{r,p \otimes q}$ and 
$H_{p \otimes q,r}$ is equal to $N$ if $r$ is equivalent to 
$pq$, and zero otherwise. In the first case, the 
intertwiners (which are injections) are called 
\emph{Clebsch-Gordan operators} (CGO), and in the second 
case they are called \emph{dual Clebsch-Gordan operators} 
(these are projectors). We shall now give 
an explicit basis of the CGO for a 
regular pair of standard representations; for the dual 
operators, see Proposition \ref{6j}.

\medskip

\noindent Consider the curve $\Gamma \subset \mathbb{CP}^2$ 
which is the zero set of the Fermat equation 
$x^N + y^N = z^N$ (with homogeneous coordinates), and define for 
any positive integer $n$ a rational function $\omega$ 
(not to be confused with the root of unity $\omega$ !) by:
$$\omega(x,y,z \vert n) = \prod_{j=1}^n \frac{y}{z-x\omega^j}
\ ,\ [x,y,z] \in \Gamma \setminus \{[1,0,\omega^j],j=1,
\ldots,n\}.$$
\noindent Moreover, set $\omega(x,y,z \vert m, n) = 
\omega(x,y,z \vert m-n)\omega^{n^2/2}$. Define also a 
periodic Kronecker symbol by $\delta(n) = 1$ if 
$n \equiv 0 \ \rm{mod} \ N$, and zero otherwise. Note that 
the function $\omega$ is periodic in its integer argument, 
with period $N$. It satisfies very nice transformation 
properties for some automorphisms of $\Gamma$, and it 
verifies summation formulae which are analogues in the root 
of unity case to well-known formulae for hypergeometric 
series. See \cite{42} for a summary, and \cite{43} and the 
references therein for details on their use in statistical 
mechanics.  

\begin{prop} \label{CG} Let $(p,q)$ be a regular pair of 
standard representations. For any non-zero complex number 
$\nu_{p,q}$, the set $\{K_\alpha(p,q), \alpha = 
0,\ldots,N-1\}$ of linear operators with components
\begin{eqnarray} K_\alpha(p,q)_{i,j}^k = \nu_{p,q} \ 
\omega^{\alpha j}\omega(t_px_q,\frac{x_p}{t_q},x_{pq}\vert i,
\alpha) \delta(i+j-k) \nonumber \end{eqnarray}
is a basis of $H_{pq, p\otimes q}$.
\end{prop}

\noindent The factor $\nu_{p,q}$ is of course inessential 
here, but it is justified below.

\medskip

\noindent Consider the following commutative diagram, where 
$(p,q,r)$ is a regular sequence of standard representations 
and the arrows denote injections of representations:
$$\begin{array}{ccc} pqr & \rightarrow & p \otimes qr \\ 
\downarrow & & \downarrow \\ pq \otimes r & \rightarrow  & p 
\otimes q \otimes r 
\end{array}$$ 
This diagram is induced, for instance, by the families of 
compositions of CGO $\{(id \otimes K_\delta(q,r)) \circ K_
\gamma(p,qr)\}_{\delta,\gamma}$ for the path of injections 
which goes to the right, and 
$\{(K_\alpha(p,q) \otimes id) \circ K_\beta(pq,r)\}_{\alpha,
\beta}$ for the path of injections which goes to the bottom. 
The \emph{$6j$-symbols} $R(p,q,r)$ of the above diagram are 
defined as the intertwiners between the two paths, i.e. it 
is a map:
\begin{eqnarray*} R(p,q,r)\ :\ H_{pqr, p\otimes qr} \otimes 
H_{qr,q\otimes r} \rightarrow H_{pq,p\otimes q} 
\otimes H_{pqr,pq\otimes r}, 
\end{eqnarray*}
which reads in the basis of CGO obtained in Proposition 
\ref{CG} as follows:
\begin{eqnarray}\label{6jdef}
K_\alpha(p,q) K_\beta(pq,r) = 
\sum_{\delta,\gamma = 0}^{N-1} R(p,q,r)_{\alpha,\beta}
^{\gamma,\delta} K_\delta(q,r) K_\gamma(p,qr).
\end{eqnarray}
These $6j$-symbols satisfy a $3$-cocycle relation which is 
easy to obtain using (\ref{6jdef}), and called the 
\emph{pentagon equation}. For a regular sequence 
$(p,q,r,s)$ of cyclic representations  of $\mathcal{W}_N$ it reads: 
\begin{eqnarray} \label{P}
R_{12}(p,q,r)R_{13}(p,qr,s)R_{23}(q,r,s) = 
R_{23}(pq,r,s)R_{12}(p,q,rs),
\end{eqnarray}
where indices denote the tensor factor on which the 
$6j$-symbols act. Next we describe the 
explicit form of these $6j$-symbols in the basis of CGO 
given in Proposition \ref{CG}.

\medskip

\noindent Introduce the functions
$$\forall x \in \mc^*,\ g(x) = \prod_{j=1}^{N-1}(1 - x
\omega^j)^{j/N},\ 
h(x) = x^{-p}\frac{g(x)}{g(1)},$$
\noindent where the branch of the $N$'th- root is chosen by 
the reality condition $g(0)=1$ (recall that $N$ is odd), and $p$ is defined by $N = 2p + 1$. Then fix the scalar $\nu_{p,q}$ in Proposition 
\ref{CG} as $\nu_{p,q} = h(\frac{x_{pq}}{t_px_q})$, and set 
$$\nu_{p,q,r} = 
h\left(\frac{x_{pq}x_{qr}}{x_{pqr}x_q}\right), \  
[x] = N^{-1}\left(\frac{1-x^N}{1-x} \right).$$ 
\noindent Note that the functions $\nu$ for the CGO are 
chosen so that, with Proposition \ref{6j}, one can write :
$$ K_\alpha(p,q)_{i,j}^k = R(p,q)_{\alpha,k}^{i,j}\ .
\nonumber$$
\noindent This is not a coincidence : one can prove that 
the $6j$-symbols and the CGO obtained above are both 
representations of the canonical element of the Heisenberg 
double of the Weyl algebra (which is a twisted quantum 
dilogarithm), acting on $H_{pqr, p\otimes qr} \otimes 
H_{qr,q\otimes r}$. This and the relation between 
(\ref{6jdef}) and (\ref{P}) are explained in \cite{16}.   

\begin{prop} \label{6j} The $6j$-symbols of the regular 
sequences of cyclic representations of $\mathcal{W}_N$, in 
the basis of the Clebsch-Gordan operators obtained in 
Proposition \ref{CG}, are described by the following 
components:
$$R(p,q,r)_{\alpha,\beta}^{\gamma,\delta} = \nu_{p,q,r} \ 
\omega^{\alpha\delta}\ \omega(x_{pqr}x_q,x_px_r,x_{pq}x_{qr}
\vert \gamma,\alpha) \ \delta(\gamma + \delta - \beta),$$
\noindent and their inverses are given by
$$\bar{R}(p,q,r)_{\gamma,\delta}^{\alpha,\beta} =  
\frac{[\frac{x_{pqr}x_{q}}{x_{pq}x_{qr}}]}{\nu_{p,q,r}}\ 
\omega^{-\alpha\delta}\ \frac{\delta(\gamma + 
\delta - \beta)}{\omega(\frac{x_{pqr}x_q}{\omega},x_px_r,
x_{pq}x_{qr}\vert \gamma,\alpha)}\ .$$
\end{prop}

\begin{remark} {\rm Restricting attention to standard 
representations $p$ where $t_p = 1$ would also give the 
above $6j$-symbols (this is how they are presented in 
\cite{1}). Moreover, one can continue them to nilpotent (i.e. non-cyclic)
representations of $\mathcal{W}_N$ with $x_p=0$, but this is of no interest for
the $K_N$ invariants.}
\end{remark}

\noindent Let $(W,L)$ be as usual. Consider a decoration 
$\mathcal{T} = (T,H,b,z,c)$ of $(W,L)$ and fix a common 
determination $d$ of the $N$'th-root for the entries of 
$\{ z(e) \}$; $x(e)$ and $y(e)=d(x(e))$ are as in \S \ref{charges}. 

\noindent Each tetrahedron $\Delta \in T$ inherits a 
decoration from $(\mathcal{T},d)$; in particular, using $b$, 
the $y(e)$'s satisfy a specific Fermat equation, and there is 
defined an index for $\Delta$, see \S \ref{branchings}. 
From now on, we denote by $x_{ij}$ (resp. $c_{ij}$) the value 
of $y$ (resp $c$) on the oriented edge (resp. the 
non-oriented edge) of $\Delta$ with vertices 
numbered $v_i$ and $v_j$ via the branching. 
Moreover, a state $\alpha : \ F(T) \to \mz/N\mz$ produces a 
number $\alpha_k$ associated to the face of $\Delta$ opposite 
to its $k$th-vertex. 

\noindent Suppose that $\Delta$ has index $-1$, and define a 
correspondence $t$ as follows: 
send $(\Delta,\mathcal{T},d,\alpha)$ to the number 
$R(\mathcal{T},d) = R(x_{01},x_{12},x_{23})_{\alpha_3,
\alpha_1}^{\alpha_2,\alpha_0}$. One verifies easily that 
under the permutation $(v_1,v_2,v_3)$ (which preserves the 
index), this number turns into the component of a matrix 
with completely different $x$-parameters. Hence we are led to 
consider the action of elementary symmetries of $\Delta$, 
considered as an abstract tetrahedron, on $t$. The 
transpositions $(v_0,v_1),(v_1,v_2),(v_2,v_3)$ are 
convenient, and they generate the whole group of symmetries 
of $\Delta$. Since a transposition of vertices changes the 
index, it is natural from a topological point of view 
(see Proposition \ref{orth} below) to associate via $t$ the 
component $\bar{R}(x_{01},x_{12},x_{23})_{\alpha_2,\alpha_0}
^{\alpha_3,\alpha_1}$ to $(\Delta,\mathcal{T},d,\alpha)$ when 
the index is $1$. Note that the charge is still inessential 
here.

\medskip

\noindent Define $N \times N$ matrices (distinguished by the 
uppered and lowered indices) by the components
$$\begin{array}{ll}
T_{m,n} = \zeta^{-1}\omega^{m^2/2}\delta(m + n), &  
\ S_{m,n} = N^{-1/2}\omega^{mn}, \\
T^{m,n} = \zeta\omega^{-m^2/2}\delta(m + n),\ \zeta \in 
\mathbb{C}^*, & \ S^{m,n} = N^{-1/2}\omega^{-mn}.
\end{array}$$
\noindent One may show that, for a suitable $\zeta$, the 
above transpositions change the $6j$-symbol associated to 
$(\Delta,\mathcal{T},d)$ via $t$ in two different ways. 
They act:\hfill\break
1) on the scalar factor,\hfill\break
2) on $\alpha$: this action ``holds'' on the faces of 
$\Delta$ (which permute under a permutation of the vertices), 
and this translates algebraically into the action of the 
matrices $S$ and $T$ on the $6j$-symbols and the 
multiplication by a power of $\omega$ depending affinely on 
$\alpha$.

\smallskip

\noindent We then have to replace $R$ and $\bar{R}$ by new 
covariant matrices. One can avoid the scalar transformation 
by multiplying  $R(p,q,r)$ and $\bar{R}(p,q,r)$ 
by $(x_{pq}x_{qr})^p$. Then, note that one can act in an affine way on the 
$\alpha$-parameters of the components of $R$ and $\bar{R}$ by:\hfill\break
- multiplying $R$ and $\bar{R}$ by a power of $\omega$ which depends 
affinely on $\alpha$,\hfill\break
- adding a scalar in the $\mathbf{Z}_N$ arguments of the $\omega$ function, 
for instance in the first one, denoted by $\gamma$ in Proposition 
\ref{6j}; then we have to add the same scalar to $\beta$, due to the 
periodic Kronecker symbol.

\noindent The following proposition describes an extension of the 
$6j$-symbols of the regular sequences of standard representations of 
$\mathcal{W}_N$ that corresponds to this description.

\begin{prop} \label{symmetry}
Given $a,c \in \mathbf{Z}_N$ and a regular sequence $(p,q,r)$ of standard 
representations, consider the two matrices called \emph{c-$6j$-symbols} 
and defined by
$$\begin{array}{l}
R(p,q,r\vert a,c)_{\alpha,\beta}^{\gamma,\delta} = (x_{pq}x_{qr})^p
\omega^{c(\gamma - \alpha) - ac/2} R(p,q,r)_{\alpha,\beta - a}^{\gamma - 
a ,\delta},\\ 
  \\
\bar{R}(p,q,r\vert a,c)_{\gamma,\delta}^{\alpha,\beta} = (x_{pq}x_{qr})
^p\omega^{c(\gamma - \alpha) + ac/2}\bar{R}(p,q,r)_{\gamma+a,\delta}^
{\alpha,\beta+a}.
\end{array}$$
\noindent We have the following relations:
$$\begin{array}{l}
\sum_{\alpha ',\gamma '=0}^{N-1} R(p,q,r\vert a,c)_{\alpha ',\beta}^
{\gamma ',\delta} T_{\gamma,\gamma '} T^{\alpha,\alpha'} =\omega^
{a/4}\bar{R}(\bar{p},pq,r\vert a,b)_{\gamma,\beta}^{\alpha,\delta},\\
\\
\sum_{\alpha ',\delta '=0}^{N-1} R(p,q,r\vert a,c)_{\alpha ',\beta}^
{\gamma,\delta'} T_{\delta,\delta '} S^{\alpha,\alpha'} =\omega^{-c/4}
\bar{R}(pq,\bar{q},qr\vert b,c)_{\beta,\delta}^{\alpha,\gamma},\\
\\
\sum_{\beta ',\delta '=0}^{N-1} R(p,q,r\vert a,c)_{\alpha,\beta '}^{\gamma,
\delta'} S_{\delta,\delta '} S^{\beta,\beta '} =\omega^{a/4}\bar{R}(p,qr,
\bar{r}\vert a,b)_{\alpha,\delta}^{\gamma,\beta},
\end{array}$$
\noindent where $b = 1/2 -a -c \in \mathbf{Z}_N$ and $\zeta = (-1)^p\omega^
{(1-N)(N-2)/24}\omega^{1/8}$. Note that $\omega^{1/4}= \omega^{p^2}$ is an 
$N$'th-root of $1$.
\end{prop}
\noindent It is now clear how to turn $t$ into a covariant correspondence 
for the action of $\mathbf{S}_4$. The charge $c$ in $\mathcal{T}$ clearly 
provides $c$-labels as in the statement. Then define:
\begin{eqnarray} \label{assoc}
 \forall \Delta \in T,\ t^{\Delta}(\mathcal{T},d,\alpha) = \left\lbrace 
\begin{array}{l}
R(x_{01},x_{12},x_{23} \vert c_{01},c_{12})_{\alpha_3,\alpha_1}^{\alpha_2,
\alpha_0}\  \rm{if \ the \ index \ of\ } \Delta \ \rm{is \ -1},\\ \\
\bar{R}(x_{01},x_{12},x_{23} \vert c_{01},c_{12})_{\alpha_2,\alpha_0}^
{\alpha_3,\alpha_1} \ \rm{otherwise}.
\end{array} \right.
\end{eqnarray}

\noindent We say that 
$$t(\mathcal{T},d) = \prod_{\Delta} t^{\Delta} (\mathcal{T},d) = 
\sum_{\alpha} \prod_{\Delta}t^{\Delta}(\mathcal{T},d,\alpha) $$
\noindent is the operator associated to $(\mathcal{T},d)$. Of course it is 
a complex number (since each face of $T$ appears twice in the set of faces 
of the abstract tetrahedra), but this definition may also be extended to 
simplicial complexes with boundary (see below). 
\begin{remark} \label{matrices} {\rm The matrices $T_{m,n}$ and $T^{m,n}$ 
(resp. $S_{m,n}$ and $S^{m,n}$) are inverse one to each other and 
$$S^4 = id,\ S^2 = \zeta'(ST)^3,\ \vert \zeta' \vert = 1.\nonumber$$ 
This uses the 
explicit calculation of the phase factor of $g(1)$, see \cite{43}. Hence 
the matrices $S$ and $T$ define a projective $N$-dimensional 
representation of $SL(2,\mathbf{Z})$, which up to an $N$'th-root of unity 
is well-known as the modular representation on the space of characters of 
minimal models in Conformal Field Theory \cite[\S 10.5]{44}.}
\end{remark}

\noindent Let us now specify the rule which determines the action of the 
matrices $S$ and $T$, when permuting vertices $v_i$ and $v_{i+1}$ of 
$\Delta$ with adjacent numberings.
\noindent  Each face $f \in F(\Delta)$ containing $v_i$ and $v_{i+1}$ 
inherits an orientation from the branching : set $\epsilon(f) = 1$ if this 
orientation is the one induced by $\Delta$ as a boundary, and 
$\epsilon(f) = -1$ otherwise. Moreover, set $\lambda(f) = 1$ if the 
numbering of the vertex of $f$ distinct from $v_i$ and $v_{i+1}$ is 
greater than $i+1$, and $\lambda(f) = -1$ if it is lesser than $i$. 
Finally, define :

$$\mu_{ab}(f) = \frac{(1 + a\lambda(f))(1 + b\epsilon(f))}{4},\ a,b = \pm. 
\nonumber$$

\noindent For example, if $\Delta$ has index $-1$ and we consider the permutation $(v_0,v_1)$, the only 
non-zero $\mu$-number associated to the face $f_2$ is $\mu_{++}$. Then, 
for any permutation $(v_i,v_{i+1})$, the matrix with components

$$s_u(f)_{m,n} = T_{m,n}\mu_{++}(f) + T^{m,n}\mu_{+-}(f) + S_{m,n}\mu_{-+}
(f) + S^{m,n}\mu_{--}(f) \nonumber$$

\noindent acts on $t^{\Delta}(\mathcal{T},d)$ by multiplication through the 
tensor factors corresponding to the $\alpha$-labels of the faces $f$ 
containing both $v_i$ and $v_{i+1}$. In particular, if we change the 
$\epsilon$-value of such an $f$, $\mu_{.,+}$ turns into $\mu_{.,-}$ etc, and the 
dual matrix acts. Then, a transposition of vertices acts on the operator 
$t(\mathcal{T},d)$ only by multiplication by a power of $\omega$ 
(see Proposition \ref{symmetry}), since for each face which is acted on 
and for the two tetrahedra glued on it the $\mu$-coefficients are opposite.

\medskip

\noindent Set 
$$Q_N(\mathcal{T},d) = \left( t(\mathcal{T},d) \right)^N = \left( 
\sum_{\alpha} \prod_{\Delta} t^{\Delta}(\mathcal{D},\alpha,d) \right)^N,$$
\noindent and consider a maximal tree $\Gamma$ in the $1$-skeleton of the 
dual cell decomposition of $T$. Denote by $\Gamma(T)$ the polyhedron 
obtained from the gluing of the abstract tetrahedra $\Delta \in T$ along 
the faces dual to the edges of $\Gamma$.
Let $\alpha_{\Gamma}$ be a state of $\Gamma(T)$, and consider the operator
$$\Gamma_N(\mathcal{T},d) = \sum_{\alpha_{\Gamma}} \prod_{\Delta} 
t^{\Delta}(\mathcal{D},\alpha,d).$$
\noindent We clearly have:
$$Q_N(\mathcal{T},d)= \left( tr(\Gamma_N(\mathcal{T},d)) \right)^N,$$
\noindent where $tr$ is the trace on $End(\mc^{\otimes e})$ and $e$ is the 
number of edges of $\Gamma$. This gives an essential 
representation-theoretic interpretation of the state sums 
$K_N(\mathcal{T},d)$ as $N$'th-powers of weighted traces, the weights being given by the 
product of the scalars $x(e)^{-2p/N},e \in T \setminus H$. It makes a bridge between 
the ``observable'' topological quantities in the invariants $K_N$ and the 
algebraic properties of the quantum dilogarithm operator on 
$\mathcal{W}_N$. In particular, it  should be a very important tool in the 
study of the asymptotic behaviour of $K_N$.

\begin{lem}\label{indbranch}
$Q_N(\mathcal{T},d)$ (whence $K_N(\mathcal{T},d)$) does not depend on the 
branching $b$ in $\mathcal{T}$.
\end{lem}

\noindent {\it Proof.} \ Any change of branching of $T$ translates on the abstract $\Delta$'s into a composition of transpositions of the vertices. We have 
seen in Proposition \ref{symmetry} that for such symmetries, the matrices 
$S$ and $T$ or their duals act on $t^{\Delta}(\mathcal{D},d)$ through the 
tensor factors corresponding to the faces $f$ containing the two permuting 
vertices; moreover, there is a power of $\omega$ appearing in factor. But 
we saw above that the matrix action is trivial on the operator 
$t(\mathcal{T},d)$. Then $Q_N(\mathcal{T},d)$ does not depend on the 
branching.

\noindent Using the above terminology, this means that a change of 
branching turns $\Gamma_N(\mathcal{T},d) \in End(\mc^{\otimes e})$, up to 
an $N$'th-root of unity, into a conjugate operator. Then its trace 
$Q_N(\mathcal{T},d)$ does not depend on $b$.

\medskip

\noindent Here is a simple consequence of the definition of the c-$6j$-symbols:

\begin{prop} \label{unitarite}
Let $(p,q,r)$ be a regular sequence of standard representations of 
$\mathcal{W}_N$. We have the following identity :
$$\bar{R}(p^*,q^*,r^*\vert a,c)_{\gamma,\delta}^{\alpha,\beta} = 
\left( R(p,q,r\vert a,c)_{-\alpha,-\beta}^{-\gamma,-\delta} \right)^*,
\nonumber$$
\noindent where $*$ denotes the complex conjugation.
\end{prop}

\noindent We now present the extension of the pentagon equation that the 
c-$6j$-symbols satisfy. It should be in the form of equation (\ref{P}), 
with some additional conditions on the charges $c$. Consider the branched 
convex polyhedron $P_0$ obtained by gluing two abstract tetrahedra along 
the face opposite to $v_1$ and $v_3$ respectively, so that 
$\Delta(v_0,v_2,v_3,v_4)$ has index $-1$ (see the left of Fig. \ref{w(e)}). Denote by $c^i$ the 
restriction of $c$ to the tetrahedron $\Delta^i$ which do not 
contain $v_i$.

\noindent We shall now prove that the operator $t(P_0,\mathcal{T},d)$ is 
equal to $t(P_1,\mathcal{T}_1,d_1)$, where $(P_1,\mathcal{T}_1,d_1)$ is 
obtained via a decorated $2 \to 3$ transit. This transit is dually 
described in the third picture of Fig. \ref{fig:slide}. The resulting 
equation is called the \emph{extended pentagon relation} (EP relation 
below).

\medskip

\noindent Let us first comment the behaviour of charges under decorated 
$2 \leftrightarrow 3$ transits. It is easy to verify that there are 
exactly four degrees of freedom on the charges of $(P_0,\mathcal{T},d))$ 
(e.g. $c_{03}^1,c_{01}^3,c_{23}^1$ and $c_{21}^3$), and there is one more 
independant degree of freedom on the charges of $(P_1,\mathcal{T}_1,d_1)$ 
(e.g. $c_{03}^4$). Some details on the geometric interpretation of the 
latter were given in  Lemma \ref{charge-transit}. Then there are five 
degrees of freedom on the charges of the EP relation. In particular, note 
that we have :
$$c_{02}^1 + c_{24}^1 + c_{40}^1 = (c_{02}^4 - c_{02}^3) + (c_{24}^0 - 
c_{24}^3) + (c_{04}^2 - c_{04}^3) = c_{13}^4 + c_{13}^0 + c_{13}^2 - 
(c_{02}^3 + c_{24}^3 + c_{40}^3),\nonumber$$
\noindent where we use the fact that opposite edges have the same charge 
in the second equality. Now this gives :
\begin{eqnarray}
c_{02}^1 + c_{24}^1 + c_{40}^1 = c_{02}^3 + c_{24}^3 + c_{40}^3 \equiv 1/2 
\pmod N \Leftrightarrow \  c_{13}^4 + c_{13}^0 + c_{13}^2 \equiv 1 \pmod N.
\label{e0}
\end{eqnarray}
\noindent Then, we see, once again, that \ref{defcharges} (1) is the 
proper setting for defining the transit of charges. Let us consider the 
following set of independent charges for the EP relation :
$$i = c_{01}^4,j = c_{01}^2,k = c_{12}^0,l = c_{23}^1,m = c_{12}^3\ .
\nonumber$$
\noindent One can easily show, using Definition \ref{defcharges} (1) 
and the charge-transit condition, that
$$l + m = c_{13}^2,l-i =  c_{23}^0,j+k= c_{02}^1, i+j =  c_{01}^3,m-k =  
c_{12}^4\ .\nonumber$$

\begin{prop} \label{EP}
Let $(p,q,r,s)$ be a regular sequence of standard representations of 
$\mathcal{W}_N$.

i) The following EP relation holds :
\begin{eqnarray}
R_{12}^4(p,q,r \vert i,m-k)R_{13}^2(p,qr,s \vert j,l+m)R_{23}^0(q,r,s 
\vert k,l-i) =  \hspace{3cm} \nonumber \\
\hspace{5cm} x_{qr}^{2p}R_{23}^1(pq,r,s \vert j+k,l)R_{12}^3(p,q,rs 
\vert i+j,m).\nonumber
\end{eqnarray}
\noindent Setting $x_p = x_{01},x_{q} =x_{12},x_r = x_{23}$ and 
$x_s = x_{34}$, the left-hand side (resp. the right-hand side) is the 
operator associated via (\ref{assoc}) to $(P_1,\mathcal{T}_1,d_1)$ (resp. 
$(P_0,\mathcal{T},d)$).
\smallskip

ii) A similar identity holds up to an $N$'th-root of unity for any other 
branchings of $P_0$ and $P_1$, even if one cannot blow down the branching of 
$P_1$ to the one of $P_0$.
\end{prop}

\noindent Remark that $ii)$ gives an alternative proof of the 
charge-invariance (3) in Proposition \ref{propo1}. Using the symmetry relations in Proposition \ref{symmetry}, one can easily derive the whole set of EP relations, for any branching of the abstract polyhedron $P_0$. 

\medskip

\noindent We give the algebraic relation that corresponds to the 
$0 \to 2$ decorated transit in the case where the branching gives the 
simplest form of it. This branching is induced from that of 
$(P_1,\mathcal{T}_1,d_1)$ above. Indeed, one can read a relation 
corresponding to a $0 \to 2$ decorated transit for $\Delta^4$ by 
comparing the identities obtained by applying first a $2 \rightarrow 3$ transit on $\Delta^0,\Delta^2 \subset P_1$, 
and then a $3 \to 2$ transit on $\Delta^0,\Delta^2$ and $\Delta^4$ (this is possible, since the final configuration is branchable). Then, 
we have:

\begin{prop} \label{orth}
Let $(p,q,r)$ be a regular sequence of standard representations of 
$\mathcal{W}_N$. The following orthogonality relation holds :
$$R(p,q,r\vert a,c)\bar{R}(p,q,r \vert -a,-c) = (x_{pq}x_{qr})^{2p} id 
\otimes id.$$
\end{prop}
\noindent Finally, consider the bubble move, which consists in 
decomposing a face of $T$ into three new ones with a common vertex inside
 and two tetrahedra glued along them. A branched relation corresponding to a bubble move is obtained by taking the partial trace over one of the 
indices in the above orthogonality relation:
\begin{prop} \label{vertex}
Let $(p,q,r)$ be a regular sequence of standard representations of 
$\mathcal{W}_N$. The following bubble relation holds:
$$tr_i \left( R(p,q,r\vert a,c)\bar{R}(p,q,r \vert -a,-c) \right) = N \ 
(x_{pq}x_{qr})^{2p} \ id,$$
where $i$ is equal to $1$ or $2$.
\end{prop}

\noindent Since any $1 \to 4$ transit may be obtained from a bubble move followed by a $2 \to 3$ transit, Propositions 
\ref{EP} and \ref{vertex} plus the symmetry relations give the whole set of relations 
corresponding to decorated $1 \to 4$ transits. This concludes this 
Appendix.

\bigskip

\noindent {\bf Aknowledgments.} The first named author is grateful for the kind
invitation of the Department of Mathematics of the University of Pisa in 
September 2000, while a part of the present work has been worked out.
The second named author profited of the kind hospitality at the University
P. Sabatier of Toulouse in December 2000. Both are indebted to Peter Zvengrowski for a careful reading of the original version of this paper.

\end{document}